\documentclass[smallextended,glov3]{svjour3}
  \def\ZR{{\mathbb R}}

\def\beq{\begin{equation}}
\def\eeq{\end{equation}}
\def\be{\begin{equation}}
\def\ee{\end{equation}}
\def\beqar{\begin{eqnarray}}
\def\eeqar{\end{eqnarray}}
\def\ber{\begin{eqnarray}}
\def\eer{\end{eqnarray}}
\def\berb{\begin{eqnarray*}}
\def\eerb{\end{eqnarray*}}

\def\SO{\mathop{\rm SO}\nolimits}
\def\so{\mathop{\rm so}\nolimits}

\def\det{\mathop{\rm det}\nolimits}

\def\tr{\mathop{\rm tr}\nolimits}

\def\norm#1.#2.{\|#1\|_{#2}}
\def\Norm#1.#2.{\big\|#1\big\|_{#2}}
\def\NOrm#1.#2.{\bigg\|#1\bigg\|_{#2}}
\def\NORm#1.#2.{\Big\|#1\Big\|_{#2}}
\def\NORM#1.#2.{\Bigg\|#1\Bigg\|_{#2}}

\def\vec#1{{\mathchoice{\mbox{\boldmath$\displaystyle#1$}}
{\mbox{\boldmath$\textstyle#1$}}
{\mbox{\boldmath$\scriptstyle#1$}}
{\mbox{\boldmath$\scriptscriptstyle#1$}}}}

\newcommand{\sym}{\mathop{\rm sym}\nolimits}

\def \0b{{\hbox{\boldmath $0$}}}

 \newcommand{\bbb}{{\hbox{\bf b}}}

\newcommand{\ubb}{{\hbox{\bf u}}}

\newcommand{\ab}{\vec{a}} 
\newcommand{\cb}{\vec{c}}

\newcommand{\db}{\vec{d}}
\newcommand{\eb}{\vec{e}} \newcommand{\fb}{\vec{f}}

 \newcommand{\nb}{\vec{n}}

\newcommand{\ub}{\vec{u}} \newcommand{\vb}{\vec{v}}
 
\newcommand{\yb}{\vec{y}}

\newcommand{\Abb}{{\bf A}} \newcommand{\Bbbb}{{\bf B}}
\newcommand{\Cbb}{{\bf C}} 
 \newcommand{\Fbb}{{\bf F}}
\newcommand{\Gbb}{{\bf G}} \newcommand{\Hbb}{{\bf H}}
\newcommand{\Ibb}{{\bf I}}

\newcommand{\Qbb}{{\bf Q}} \newcommand{\Rbb}{{\bf R}}

\newcommand{\Gb}{\vec{G}}

\newcommand{\Ub}{\vec{U}}

\def \Thetab{\vec{\Theta}}

\def \Thetab{\vec{\Theta}}

\newcommand{\at}{{\tilde a}}

\newcounter{primjer}[section]
\newcounter{tvrdnja}[section]
\newcounter{zadatak}[section]
\usepackage{esint}
\usepackage{amssymb}
\usepackage{graphicx}
\usepackage{enumerate}
\usepackage{verbatim}
\newcommand{\dist}{\mathop {\mbox{\rm dist}}\nolimits}

\newcommand{\homeg}{\hat{\Omega}^h}
\newcommand{\wirbb}{\widetilde{\Rbb}}
\newcommand{\ide}{\textrm{ {\it id}}}

\begin{document}

\title{Shallow shell models  by $\Gamma$-convergence
}




\author{
        Igor Vel\v{c}i\'{c} 
}
\institute{ Igor Vel\v{c}i\'{c} \at Faculty of Electrical Engineering and Computer Science, University of Zagreb, Unska 3, 10000 Zagreb, Croatia \\Tel: +385-1-6129965 \\Fax:+385-1-6170007 \\ \email{igor.velcic@fer.hr}}



\maketitle

\begin{abstract}
In this paper we derive, by means of $\Gamma$-convergence, the
shallow shell models starting from  non linear three
dimensional elasticity. We use the approach analogous to the
one for shells and plates. We start from the minimization
formulation of the general three dimensional elastic body which
is subjected to normal volume forces and free boundary
conditions  and do not presuppose any constitutional behavior.
To derive the model we need to propose how is the order of
magnitudes of the external loads related to the thickness of
the body $h$ as well as the order of the "geometry" of the
shallow shell. We analyze the situation when the external
normal forces are of order $h^\alpha$, where $\alpha>2$. For
$\alpha=3$ we obtain the Marguerre-von K\'{a}rm\'{a}n model and
for $\alpha>3$ the linearized Marguerre-von K\'{a}rm\'{a}n
model. For $\alpha \in (2,3)$ we are able to obtain only the
lower bound for the $\Gamma$-limit. This is analogous to the
recent results for the ordinary shell models.
\end{abstract}

\section{Introduction}

The study of thin structures is the subject of numerous works
in the theory of elasticity.  Many authors have proposed
two-dimensional shell and plate models and we come to the
problem of their justification. There is a vast literature on
the subject of plates and shells and some also devoted to
shallow shells (see \cite{Ciarlet0,Vorovich}). The expression
shallow shell means that the curvature of the mean surface is
also small with respect to the sizes of the mean surface.

The derivation and justification of the lower dimensional
models, equilibrium and dynamic, of rods, curved rods, plates,
shells, shallow shells in linearized elasticity by using formal
asymptotic expansion is well established (see
\cite{Ciarlet0,Ciarlet1} and the references therein). In all
these approaches one starts from the equations of
three-dimensional linearized elasticity and then via formal
asymptotic expansion justify the lower dimensional models. One
can also obtain the convergence results. Complete asymptotic
(with higher order terms) for linear shallow shells and the
influence of boundary conditions (boundary layer) is discussed
in \cite{Andreoiu}.

Formal asymptotic expansion is also applied to derive non
linear models of plates and shells (see
\cite{Ciarlet0,Ciarlet1} and the references therein), starting
from three-dimensional isotropic elasticity (usually
Saint-Venant-Kirchoff material). Hierarchy of the models is
obtained, depending on the boundary conditions  and the order
of the external loads related to the thickness of the body $h$
(see also \cite{Fox} for plates). Formal asymptotic expansion
is also applied to derive non linear shallow shell models (see
\cite{Ciarlet0,Andreoiu2,Hamdouni}). It turns out that the
asymptotic analysis of the equations of shallow shells closely
follows that of a plate. If we assumed that the order of the
normal external loads behaves like $h^3$, we would obtain
Marguerre-von K\'{a}rm\'{a}n equations. Influence of different
lateral boundary  conditions on the model is also discussed
(see \cite{Ciarlet2,Gratie}).

However, formal asymptotic expansion does not provide us a
convergence result. The first convergence result, in deriving
lower dimensional models from three-dimensional non linear
elasticity, is obtained applying $\Gamma$-convergence, very
powerful tool introduced by Degiorgi (see
\cite{Braides,dalmaso}). Using $\Gamma$-convergence, membrane
plate and membrane shell models are obtained (see
\cite{Ledret1,LeDret2}). It is assumed that the external loads
are of order $h^0$. The obtained models are different from
those ones obtained by the formal asymptotic expansion, at
least for some specific  deformations (compression).

Recently, hierarchy of models of plates and shells is obtained
via $\Gamma$- convergence (see
\cite{Muller0,Muller3,Muller4,Muller6,Lewicka1}). Influence of
the boundary conditions and the order and  the type of the
external loads is largely discussed (see
\cite{Muller3,Muller5}). Let us mention that
$\Gamma$-convergence results provide us the convergence of the
global minimizers of the energy functional. Recently,
compensated compactness arguments are used to obtain the
convergence of the stationary points of the energy functional
(see \cite{Mora,Muller7}).

Here we apply the tools developed for plates and ordinary
shells to obtain shallow shell models by $\Gamma$-convergence.
It is assumed that we have free boundary conditions and that
the external loads are normal forces with order $\alpha$
greater than $2$. The main result is given in Theorem
\ref{najglavnijiteorem}. For the normal volume forces of order
$\alpha=3$ we obtain Marguerre-von K\'{a}rm\'{a}n  equations
and for the order $\alpha>3$ we obtain linearized Marguerre-von
K\'{a}rm\'{a}n  equations. Thus we have also justified
linearized theory from three-dimensional non linear elasticity
postulating only the smallness of the forces. We do not
presuppose any constitutional behavior of the material.

In the situations when we have specific geometry of the shallow
shell (developable surface) we are also able to obtain the
model for $\alpha \in (2,3)$. Otherwise,  greater influence of
the geometry of the shallow shell is expected (see Remark
\ref{napomenaoost}).

Throughout the paper  $\bar{A}$ or $\{A\}^-$  denotes the
closure of the set. By a domain we call a bounded open set with
Lipschitz boundary. $\Ibb$  denotes the identity matrix, by
$\SO(3)$ we denote the rotations in $\ZR^3$ and by $\so(3)$ the
set of antisymmetric matrices $3 \times 3$. $x'$ stands for
$(x_1,x_2)$. $\eb_1,\eb_2,\eb_3$ are the vectors of the
canonical base in $\ZR^3$. By $\ide$ we denote the identity
mapping $\ide (x)=x$. By $\nabla_h$ we denote
$\nabla_h=\nabla_{\eb_1,\eb_2}+\frac{1}{h} \nabla_{\eb_3}$.
$\rightarrow$ denotes the strong convergence and
$\rightharpoonup$ the weak convergence. We suppose that the
Greek indices $\alpha,\beta$ take the values in the set
$\{1,2\}$ while the Latin indices $i,j$ take the values in the
set $\{1,2,3\}$.

\section{Setting up the problem}
\setcounter{equation}{0}
We consider a three-dimensional elastic shell occupying in its reference configuration the set $\{\hat{\Omega}^h\}^-$, where $\hat{\Omega}^h=\Thetab^h(\Omega^h), \ \Omega^h=\omega \times (-h,h),\  \omega$ is a domain  in $\ZR^2$ and the mapping $\Thetab^h: \{\hat{\Omega}^h\}^- \to \ZR^3$ is given by
$$ \Thetab^h(x^h)=(x_1,x_2, \theta^h(x_1,x_2))+x_3^h \ab_3^h(x_1,x_2)$$
for all $x^h=(x',x_3^h) \in \bar{\Omega}^h$, where $\ab_3^h$ is
a unit normal vector to the middle surface
$\Thetab^h(\overline{\omega})$ of the shell. We assume that
$\theta^h=f^h\theta$, where $\theta \in C^2(\bar{\omega})$,
$\lim_{h \to 0} f^h=0$, $f^h>0$. Using that assumption we can
conclude that at each point of the surface $\bar{\omega}$ the
vector $\ab_3^h$ is given by
$$ \ab_3^h=(\alpha^h)^{-1/2}(-f^h\partial_1 \theta, -f^h \partial_2 \theta, 1),$$
where  $$\alpha^h=(f^h)^2|\partial_1
\theta|^2+(f^h)^2|\partial_2 \theta|^2+1.$$ By inverse function
theorem it can be easily seen that for $h\leq h_0$ small enough
$\Thetab^h$ is a $C^1$ diffeomorphism  (the global injectivity
can be proved by adapted compactness argument, see  \cite[Thm
3.1-1]{Ciarlet1} for the ordinary shell). The following theorem
is easy to prove and is a direct consequence of Theorem 3.3-1.,
page 219, 220 in \cite{Ciarlet0}.
\begin{theorem} \label{izcea}
Let the function $\theta^h$ be such that
$$\theta^h (x_1,x_2)=f^h \theta (x_1,x_2), \mbox{ {\rm for all} } (x_1,x_2) \in \bar{\omega}, $$
where $\theta \in C^2(\bar{\omega})$ is independent of $h$.
Then there exists $h_0=h_0(\theta)>0$ such that the Jacobian
matrix $\nabla \Thetab^h(x^h)$ is invertible for all $x^h \in
\bar{\Omega}^h$ and all $h \leq h_0$. Also there exists $C>0$
such that  for $h \leq h_0$ we have
\begin{equation} \label{determinanta} \det \nabla \Thetab^h=1+(f^h)^2 \delta^h (x^h), \end{equation}
and
\begin{equation} \label{nablaje}
\nabla \Thetab^h(x^h)=\Ibb-f^h \Cbb(x')+\max\{(f^h)^2, h(f^h) \} \Rbb_1^h(x^h),
\end{equation}
\begin{equation} \label{inverz}
 (\nabla \Thetab^h(x^h))^{-1}=\Ibb+f^h \Cbb(x') +\max\{(f^h)^2, h(f^h) \} \Rbb_2^h(x^h),
\end{equation}
\begin{equation} \label{glupaocjena}
\| (\nabla \Thetab^h) -\Ibb \|_{L^{\infty}(\Omega^h;\ZR^{3 \times 3})},\| (\nabla \Thetab^h)^{-1} -\Ibb \|_{L^{\infty}(\Omega^h;\ZR^{3 \times 3})}<Cf^h,
\end{equation}
\begin{equation} \label{kristina}
\left\|\frac{1}{hf^h} \frac{\nabla \Thetab^h(x',x_3^h+hs)-\nabla \Thetab^h(x',x_3^h)}{s}\right\|_{L^\infty} \leq C,
\end{equation}
where
\begin{equation} \label{defA}
\Cbb(x')=\left( \begin{array}{ccc} 0 & 0 & \partial_1 \theta(x') \\ 0 & 0 & \partial_2 \theta(x') \\ -\partial_1 \theta(x') & -\partial_2 \theta(x') & 0 \end{array} \right)
\end{equation}
and $\delta^h: \bar{\Omega}^h \to \ZR, \
\Rbb_k^h:\bar{\Omega}^h \to \ZR^{3 \times 3},\ k=1,2$ are
functions which satisfy
$$ \sup_{0<h \leq h_0} \max_{x^h \in \bar{\Omega}^h} |\delta^h(x^h)| \leq C_0, \sup_{0<h \leq h_0} \max_{i,j} \max_{x^h \in \bar{\Omega}^h} |\Rbb_{k,ij}^h(x^h)| \leq C_0, \ k=1,2,$$
for some constant $C_0>0$.
\end{theorem}
\begin{prooof}
It is easy to see (see also \cite[p.220]{Ciarlet0})
\begin{eqnarray} \nonumber
(\nabla \Thetab^h (x^h))_{11} &=& 1-\frac{f^h}{2} x_3^h (\alpha^h (x'))^{-3/2}(2 \alpha^h (x') \partial_{11} \theta(x')\\ \label{rel1} & &-\partial_1 \alpha^h (x') \partial_1 \theta (x')), \\
(\nabla \Thetab^h (x^h) )_{12} &=& -\frac{f^h}{2} x_3^h(\alpha^h(x')^{-3/2} (2 \alpha^h (x') \partial_{12} \theta-\partial_2 \alpha^h (x') \partial_1 \theta(x')), \\
(\nabla \Thetab^h (x^h))_{13} &=& -f^h (\alpha^h(x'))^{-1/2} \partial_1 \theta (x'), \\ \nonumber
(\nabla \Thetab^h (x^h))_{21} &=& -\frac{f^h}{2} x_3^h (\alpha^h (x'))^{-3/2}(2 \alpha^h (x') \partial_{12} \theta(x')\\ & & -\partial_1 \alpha^h (x') \partial_2 \theta (x')), \\ \nonumber
(\nabla \Thetab^h (x^h))_{22} &=& 1-\frac{f^h}{2} x_3^h (\alpha^h (x'))^{-3/2}(2 \alpha^h (x') \partial_{22} \theta(x')\\ & &-\partial_2 \alpha^h (x') \partial_2 \theta (x')), \\
(\nabla \Thetab^h (x^h))_{23} &=& -f^h (\alpha^h(x'))^{-1/2} \partial_2 \theta (x'), \\
(\nabla \Thetab^h (x^h))_{31} &=& f^h \partial_1 \theta (x') -\frac{x_3^h}{2} (\alpha^h(x'))^{-3/2} \partial_1 \alpha^h(x'), \\
(\nabla \Thetab^h (x^h))_{32} &=& f^h \partial_2 \theta (x') -\frac{x_3^h}{2} (\alpha^h(x'))^{-3/2} \partial_2 \alpha^h(x'), \\ \label{rel9}
(\nabla \Thetab^h (x^h))_{33} &=& (\alpha^h (x'))^{-1/2}.
\end{eqnarray}
Everything is already proved in \cite{Ciarlet0}, except the relation (\ref{kristina}), which is an easy consequence of the relations (\ref{rel1})-(\ref{rel9}).
\end{prooof}

The starting point of our analysis is the minimization problem
for the shallow shell. The strain energy of the shallow shell
is given by
$$
K^{h}(\yb)= \int_{\hat{\Omega}^h} W(\nabla \yb(x))dx,
$$
where $W: \mathbb{M}^{3 \times 3} \to[0,+\infty]$ is the stored
energy density function.  $W$ is Borel measurable and, as in
\cite{Muller0,Muller3,Muller4}, supposed to satisfy
\begin{enumerate}[i)]
\item $W$ is of class $C^2$ in a neighborhood of $\SO(3)$;
\item $W$ is frame-indifferent, i.e., $W(\Fbb)=W(\Rbb\Fbb)$ for every $\Fbb \in \ZR^{3 \times 3}$ and $\Rbb \in \SO(3)$;
\item $W(\Fbb) \geq C_W \dist^2(\Fbb,\SO(3))$, for some $C_{W}>0$ and all $\Fbb \in \ZR^{3 \times 3}$, $\ W(\Fbb)=0$ if $\Fbb \in \SO(3)$.
\end{enumerate}
By $Q_3:\ZR^{3 \times 3} \to \ZR$ we denote the quadratic form $Q_3(\Fbb)= D^2 W(\Ibb)(\Fbb,\Fbb)$ and by $Q_2: \ZR^{2 \times 2} \to \ZR$ the quadratic form,
\begin{equation}
Q_2 (\Gbb)= \min_{\ab \in \ZR^3} Q_3 (\Gbb+\ab \otimes \eb_3+\eb_3 \otimes \ab),
\end{equation}
obtained by minimizing over the stretches in the $x_3$ directions. Using ii) and iii) we conclude that both forms are positive semi-definite (and hence convex), equal to zero on antisymmetric matrices and depend only on the symmetric part of the variable matrix, i.e. we have
\begin{equation}
Q_3(\Gbb)=Q_3(\sym \Gbb), \quad Q_2 (\Gbb)= Q_2 (\sym \Gbb).
\end{equation}
Also, from ii) and iii), we can conclude that both forms are positive definite (and hence strictly convex) on symmetric matrices.
For the special case of isotropic elasticity we have
\begin{eqnarray}
\nonumber Q_3 (\Fbb) &=& 2 \mu | \frac{\Fbb+\Fbb^T}{2}|^2 + \lambda (\tr \Fbb)^2, \\
Q_2 (\Gbb) &=& 2 \mu | \frac{\Gbb+\Gbb^T}{2} |^2+ \frac{2 \mu \lambda}{2 \mu + \lambda} (\tr \Gbb)^2.
\end{eqnarray}

We suppose that the external loads are dead normal volume loads and thus we have that the total energy functional is given by
$$J^h (\yb)=K^{h}(\yb)-\int_{\hat{\Omega}^h} \fb_3^h  \yb_3,$$
where  $\fb_3^h \in L^2(\hat{\Omega}^h)$. We  suppose that the
body is free at the boundary and the total energy functional is
defined on the space $W^{1,2} (\Omega;\ZR^3)$. Since the volume
of physical domain decreases with the order 1 as $h\to 0$ it is
natural to look for the $\Gamma$-limit of the sequence
$\frac{1}{h}J^h$. In fact, since the model crucially depends on
the assumption how is the order of magnitudes of external loads
related to the thickness of the body (see
\cite{Muller0,Muller3,Muller4}), we shall look for the
$\Gamma$-limit of $\frac{1}{E^h} \frac{1}{h} J^h$. We shall
analyze the situations when $h^{-2}E^h \to 0$. For the applied
forces we suppose
\begin{equation} \label{sile}
\frac{1}{h^\alpha}  \fb_3^h  \rightharpoonup \fb_3 \ \textrm{in} \ L^2(\omega),
\end{equation}
where $\alpha>2$. The following questions have to be answered:
For given order of external loads what is the order of the
energy functional such that we have non trivial $\Gamma$-limit?
How does the limit functional look like? The main result is
given by Theorem \ref{najglavnijiteorem}, where is answered to
these questions for special cases of $\alpha$. We take the
special form of applied forces where the components of the
force in $\eb_1$ and $\eb_2$ direction vanish i.e. we suppose
$\fb_1^h=\fb_2^h=0$. For the analysis of the different
situation see Remark \ref{inplane1} and Remark \ref{inplane2}.
 Also we assume
 \begin{equation} \label{uvjet}
 \int_{\hat{\Omega}^h} \fb_3^h dx=0.  
 \end{equation}
This means that  the force on the body is equal to $0$ and is a
necessary condition, since it avoids the absence of a lower
bound of the total energy functional arising from the trivial
invariance $ \yb \to \yb+\textrm{const}$. Firstly, we analyze
how the order of the strain energy affects the limit
displacement. Let us denote $I^h:=\frac{1}{h} K^h$.
\section{$\Gamma$-convergence}
\setcounter{equation}{0}
We shall need the following theorem
which can be found in \cite{Muller0}.
\begin{theorem}[on geometric rigidity]\label{tgr}
Let $U \subset \ZR^m$ be a bounded Lipschitz domain, $m \geq 2$. Then there exists a constant $C(U)$ with the following property: for every $\vb\in W^{1,2}(U;\ZR^m)$ there is associated rotation $\Rbb\in \SO(m)$ such that
\begin{equation}\label{gr}
\|\nabla \vb-\Rbb\|_{L^2(U)} \leq C(U)\|\dist(\nabla \vb,\SO(m)\|_{L^2(U)}.
\end{equation}
The constant $C(U)$ can be chosen uniformly for a family of domains which are Bilipschitz equivalent with controlled Lipschitz constants.
The constant $C(U)$ is invariant under dilatations.
\end{theorem}
In the sequel we suppose $h_0 \geq \frac{1}{2}$ (see Theorem
\ref{izcea}). If this was not the case, what follows could be
easily adapted. Let us by $P^h: \Omega \to \Omega^{h} $ denote
the map $P^h(x',x_3)=(x', 2h x_3)$. In the same way as in
\cite[Theorem 10]{Muller3} (see also \cite[Lemma
8.1]{Lewicka1}) we can prove the following theorem. For the
adaption we only need Theorem \ref{gr} and the facts that
$C(U)$ can be chosen uniformly for Bilipschitz equivalent
domains and  that the norms $\|\nabla \Thetab^h\|$, $\|(\nabla
\Thetab^h)^{-1} \|$ are uniformly bounded on $\Omega^h$ for $h
\leq \frac{1}{2}$.
\begin{theorem} \label{prepisano}
Let $\omega \subset \ZR^2$ be a domain. Let $\Thetab^h$ be as above and let  $h \leq \frac{1}{2}$. Let $\yb^h \in W^{1,2}(\hat{\Omega}^h;\ZR^3)$
and
$$E=\frac{1}{h} \int_{\homeg}\dist^2 (\nabla \yb^h,\SO(3)) dx.$$
Then there exist maps $\Rbb^h:\omega \to \SO(3)$ and
$\widetilde{\Rbb}^h:\omega \to \ZR^{3 \times 3}$, with
$|\wirbb^h| \leq C$, $\wirbb^h \in W^{1,2}(\omega, \ZR^{3
\times 3})$ such that
\begin{equation}
\| (\nabla \yb^h) \circ \Thetab^h\circ P^h -\Rbb^h \|_{L^2(\Omega)} \leq C\sqrt{E}, \quad \|\Rbb^h-\wirbb^h \|_{L^2(\omega)} \leq C\sqrt{E},
\end{equation}
\begin{equation} \label{ocjena1}
\| \nabla' \wirbb^h \|_{L^2(\omega)} \leq Ch^{-1} \sqrt{E}, \quad \|\Rbb^h-\wirbb^h\|_{L^{\infty}(\omega)} \leq Ch^{-1} \sqrt{E}.
\end{equation}
Moreover there exist a constant rotation $\bar{\Qbb}^h \in \SO(3)$ such that
\begin{equation} \label{gornje}
\| (\nabla \yb^h) \circ \Thetab^h \circ P^h- \bar{\Qbb}^h\|_{L^2(\Omega)} \leq Ch^{-1} \sqrt{E},
\end{equation}
and
\begin{equation} \label{pocjena}
\| \Rbb-\bar{\Qbb}^h \|_{L^p(\omega)} \leq C_p h^{-1} \sqrt{E}, \quad \forall p<\infty.
\end{equation}
Here all constants depend only on S (and on $p$ where indicated) and $\Omega:=\Omega^{1/2}$.
\end{theorem}
\begin{remark} \label{remmull}
By $(\nabla \yb^h) \circ \Thetab^h \circ P^h$ we have denoted $\nabla \yb^h$ evaluated at the point $\Thetab^h(P^h(x))$.
The construction of $\Rbb^h$ and $\wirbb^h$ is given in \cite{Muller3}. Since $\SO(3)$ is a smooth manifold there exists a tubular neighborhood  $\mathcal{U}$ of $\SO(3)$ such that the nearest point projection $\pi:\mathcal{U}\to \SO(3)$ is smooth.
If $E^h \leq \delta h^2$ then we always have $\wirbb^h(x') \in \mathcal{U}$. Hence the map $\Rbb_p^h:S \to \SO(3)$, given by $\Rbb_p^h(x')=\pi(\wirbb^h(x'))$,  is well defined and in Theorem \ref{prepisano} the map $\wirbb$ can be replaced by the map $\Rbb_p^h$.
This is already  noted in \cite{Muller3}.
\end{remark}
\begin{remark}
Since $\Thetab$ is Bilipschitz map, it can easily be seen that
the map $\yb \to \yb \circ \Thetab^h$ is an isomorphism between
the spaces $W^{1,2}(\Omega^h; \ZR^m)$ and
$W^{1,2}(\hat{\Omega}^h; \ZR^m)$ (see e.g. \cite{Adams}).
\end{remark}
To prove $\Gamma$-convergence result we need to prove the lower and the upper bound.
\subsection{Lower bound}
We need the following version of Korn's inequality which is proved in a standard way by contradiction.
\begin{lemma} \label{Kornova}
Let $\Omega \subset \ZR^2$ be a Lipschitz domain. Then there exists $C(\Omega)>0$ such that for an arbitrary $\ub \in W^{1,2}(\Omega; \ZR^2)$ we have
\begin{equation} \label{Korn}
\| \ub \|_{W^{1,2}(\Omega;\ZR^2)} \leq C(\Omega) (\|\sym \nabla \ub \|_{L^2(\Omega;\ZR^2)}+|\int_{\Omega} \ub dx|+|\int_{\Omega} (\partial_2 \ub_{1}-\partial_1 \ub_{2})dx|).
\end{equation}
\end{lemma}
\begin{lemma} \label{osnovna}
Let $\yb^h \in W^{1,2}(\hat{\Omega}^h;\ZR^3)$ be such that
\begin{eqnarray}
\frac{1}{h} \int_{\homeg}\dist^2 (\nabla \yb^h,\SO(3)) dx & \leq & C E^h, \\
\lim_{h\to 0} h^{-2} E^h &=&0.
\end{eqnarray}
Let us also take $f^h=\max\{h,h^{-1} \sqrt{E^h}\}$.
Then there exists maps $\Rbb^h$  $\in$ \\ $W^{1,2}(\omega, \SO(3))$ and constants $\bar{\Rbb}^h \in \SO(3)$, $\cb^h \in \ZR^3$  such that
$$ \widetilde{\yb}^h:=(\bar{\Rbb}^h)^T \yb^h- \cb^h$$
and the in-plane and the out-of-plane displacements
\begin{eqnarray} \nonumber \Ub^h (x')&:=&
\int_{-1/2}^{1/2} (\left( \begin{array}{c} \widetilde{\yb}^h_1 \circ \Thetab^h \circ P^h \\  \widetilde{\yb}^h_2 \circ \Thetab^h \circ P^h \end{array} \right)(.,x_3)-x') dx_3, \\   \label{definicijauv} V^h(x')&:=&
 \int_{-1/2}^{1/2} ((\widetilde{\yb}_3^h \circ \Thetab^h\circ P^h)(.,x_3)-f^h \theta(.)) dx_3 \end{eqnarray}
satisfy
\begin{equation} \label{sranje1}
\| (\nabla \widetilde{\yb}^h) \circ \Thetab^h \circ P^h-\Rbb^h \|_{L^2(\Omega)} \leq C \sqrt{E^h},
\end{equation}
\begin{equation}\label{sranje2}
\| \Rbb^h-\Ibb\|_{L^p(\omega)} \leq C_p h^{-1}\sqrt{E^h} \quad \forall p<\infty, \quad \| \nabla' \Rbb^h \|_{L^2(\omega)} \leq Ch^{-1} \sqrt{E^h}.
\end{equation}
Moreover every subsequence (not relabeled)  has its subsequence (also not relabeled) such that
\begin{equation}\label{konvergencija}
v^h:=\frac{h}{\sqrt{E^h}} V^h \to v  \quad \textrm{in } W^{1,2}(\omega), \quad v \in W^{2,2}(\omega),
\end{equation}
\begin{equation} \label{konvergencija2}
\ub^h:= \min\Big( \frac{h^2}{E^h}, \frac{1}{\sqrt{E^h}} \Big) \Ub^h \rightharpoonup \ub \quad \textrm{in } W^{1,2}(\omega;\ZR^2),
\end{equation}
\begin{equation} \label{acc}
\frac{h}{\sqrt{E^h}}(\Rbb^h-\Ibb) \to \Abb \quad \textrm{in } L^q(\omega;\ZR^{3 \times 3}), \quad \forall q< \infty,
\end{equation}
\begin{equation} \label{graa}
\frac{h}{\sqrt{E^h}} ((\nabla \widetilde{\yb}^h)\circ \Thetab^h \circ P^h-\Ibb) \to \Abb \quad \textrm{in } L^2(\Omega;\ZR^{3 \times 3})
\end{equation}
\begin{equation}
\partial_3 \Abb=0, \quad \Abb \in W^{1,2} (\omega; \ZR^{3 \times 3}),
\end{equation}
\begin{equation}
\Abb=\eb_3 \otimes \nabla' v -\nabla' v \otimes \eb_3,
\end{equation}
\begin{equation} \label{simetricno}
\frac{h^2}{E^h} \sym (\Rbb^h-\Ibb) \to \frac{\Abb^2}{2} \quad \textrm{in } L^2(\Omega;\ZR^{3 \times 3}).
\end{equation}
\end{lemma}
\begin{prooof}
We shall follow the proof of Lemma 13 in \cite{Muller3}.
Estimates (\ref{sranje1}) and (\ref{sranje2}) follow immediately from Theorem \ref{prepisano} and Remark \ref{remmull} since one can choose
$\bar{\Rbb}^h$ so that (\ref{gornje}) holds with $\bar{\Qbb}=\Ibb$.
Using (\ref{sranje1}) and (\ref{sranje2}) we conclude that
\begin{equation} \label{sranje3}
\| (\nabla \widetilde{\yb}^h) \circ \Thetab^h \circ P^h-\Ibb \|_{L^2(\Omega)} \leq C h^{-1}\sqrt{E^h}.
\end{equation}
For adapting the proof to the proof of Lemma 13 in \cite{Muller3}  it is essential to see
\begin{eqnarray} \nonumber
(\nabla \widetilde{\yb}^h) \circ \Thetab^h \circ P^h&=&(\nabla(\widetilde{\yb}^h \circ \Thetab^h)\circ P^h) ((\nabla \Thetab^h)^{-1} \circ P^h)\\ \label{kljucnozadokaz1}&=&\nabla_h (\widetilde{\yb}^h \circ \Thetab^h \circ P^h) ((\nabla \Thetab^h)^{-1} \circ P^h).
\end{eqnarray}
From (\ref{kljucnozadokaz1}) it follows
\begin{equation} \label{kljucnozadokaz2}
((\nabla \widetilde{\yb}^h) \circ \Thetab^h \circ P^h)((\nabla \Thetab^h) \circ P^h)=\nabla_h (\widetilde{\yb}^h \circ \Thetab^h \circ P^h).
\end{equation}
Using the fact that $\| \nabla \Thetab^h\|_{L^{\infty}(\Omega^h)}$ is bounded, (\ref{sranje3}), (\ref{kljucnozadokaz2}) we
conclude that
\begin{equation} \label{zazazaKorna}
\| \nabla_h (\widetilde{\yb}^h \circ \Thetab^h \circ P^h-\Thetab^h \circ P^h) \|_{L^2(\Omega)} \leq C h^{-1}\sqrt{E^h}
\end{equation}
 By applying additional constant in-plane rotation of order $h^{-1} \sqrt{E^h}$ to $\widetilde{\yb}^h$ and $\Rbb^h$  we may assume in addition to (\ref{sranje1}) and (\ref{sranje2}) that
\begin{equation} \label{josmalo}
\int_{\Omega} ((\partial_2 \widetilde{\yb}^h_{1})\circ \Thetab^h \circ P^h-(\partial_1 \widetilde{\yb}^h_{2})\circ \Thetab^h \circ P^h) dx=0.
\end{equation}
By choosing $\cb^h$ suitably we may also assume that
\begin{equation} \label{nacinzac}
\int_{\Omega} (\widetilde{\yb}^h\circ \Thetab^h \circ P^h-\Thetab^h \circ P^h)dx=0
\end{equation}
From (\ref{kljucnozadokaz1}) it can be easily seen that the following identity is valid
\begin{eqnarray} \nonumber
(\nabla \widetilde{\yb}^h) \circ \Thetab^h \circ P^h-\Ibb&=&(\nabla_h (\widetilde{\yb}^h \circ \Thetab^h \circ P^h-\Thetab^h \circ P^h))((\nabla \Thetab^h)^{-1} \circ P^h-\Ibb)\\ \label{zaKorna} & &+(\nabla_h (\widetilde{\yb}^h \circ \Thetab^h \circ P^h-\Thetab^h \circ P^h)).
\end{eqnarray}
Using (\ref{glupaocjena}), (\ref{zazazaKorna}),  (\ref{josmalo}), (\ref{zaKorna}) we conclude that
\begin{equation} \label{zazaKorna}
|\int_{\Omega} (\partial_2 (\widetilde{\yb}_1^h\circ \Thetab^h \circ P^h-\Thetab_1^h \circ P^h)-\partial_1 ( \widetilde{\yb}_2^h\circ \Thetab^h \circ P^h-\Thetab_2^h \circ P^h)) dx| <Cf^h h^{-1} \sqrt{E^h}.
\end{equation}
Let us define $\Abb^h=(h/\sqrt{E^h})(\Rbb^h-\Ibb)$. From  (\ref{sranje2}) we get for a subsequence
$$ \Abb^h \rightharpoonup \Abb \quad \textrm{in } W^{1,2}(\omega;\ZR^{3 \times 3}) $$
Using the Sobolev embedding we deduce (\ref{acc}).  Using (\ref{sranje1}) we deduce (\ref{graa}). Since $(\Rbb^h)^T \Rbb^h=\Ibb$ we have
$\Abb^h+(\Abb^h)^T=-(\sqrt{E^h}/h)(\Abb^h)^T \Abb^h$. Hence $\Abb+\Abb^T=0$ and after multiplication by $\frac{h}{\sqrt{E^h}}$ we obtain (\ref{simetricno}) from the strong convergence of $\Abb^h$.
From (\ref{glupaocjena}), (\ref{zazazaKorna}) and (\ref{zaKorna}) we conclude that
\begin{equation}\label{prilagodbba}
\frac{h}{\sqrt{E^h}} \nabla_h (\widetilde{\yb}^h \circ \Thetab^h \circ P^h-\Thetab^h \circ P^h) \to \Abb \quad \textrm{in } L^2(\Omega;\ZR^{3 \times 3}).
\end{equation}
From  (\ref{nacinzac}), (\ref{prilagodbba}) and the Poincare inequality we conclude the convergence in (\ref{konvergencija}).  Moreover we have $\partial_i v=\Abb_{3i}$ for $i=1,2$. Hence $v \in W^{2,2}$ since $\Abb \in W^{1,2}$. Since $\Abb$ is skew-symmetric we immediately have $\Abb_{31}=\partial_1 v, \Abb_{32}=\partial_2 v$.
If there exists constant $C>0$ such that $h^4 \leq C E^h, \forall h \leq h_0$, we could multiply (\ref{sranje1}) with $\frac{h^2}{E^h}$ to conclude
for some $C>0$
\begin{equation}
\| \frac{h^2}{E^h}\sym ((\nabla \yb^h) \circ \Thetab^h\circ P^h-\Ibb) - \frac{h^2}{E^h}\sym (\Rbb^h-\Ibb) \|_{L^2(\Omega)} \leq C.
\end{equation}
Using (\ref{simetricno}) we obtain that there exists $C>0$ such
that
\begin{equation} \label{omedjenost}
\| \frac{h^2}{E^h}\sym ((\nabla \yb^h) \circ \Thetab^h\circ P^h-\Ibb)\|_{L^2(\Omega)} \leq C.
\end{equation}
Using the identity (\ref{zaKorna}) we conclude that there exists $C>0$ such that
\begin{equation} \label{nakkkkk}
\| \frac{h^2}{E^h}\sym (\nabla_h (\yb^h \circ \Thetab^h\circ P^h- \Thetab^h\circ P^h)\|_{L^2(\Omega)} \leq C.
\end{equation}
From that we conclude that $\| \frac{h^2}{E^h} \sym \nabla' \ub^h \|_{L^2(\omega)}$ is bounded.
Using Lemma \ref{Kornova}, (\ref{nacinzac}) and  (\ref{zazaKorna}) we have the convergence (\ref{konvergencija2}).
In the situation where there exists subsequence such that $\lim_{h \to 0} h^{-4} E^h=0$ we divide (\ref{sranje1}) with $\frac{1}{\sqrt{E^h}}$ to conclude
\begin{equation}
\| \frac{1}{\sqrt{E^h}}\sym ((\nabla \yb^h) \circ \Thetab^h\circ P^h-\Ibb) - \frac{1}{\sqrt{E^h}}\sym (\Rbb^h-\Ibb) \|_{L^2(\Omega)} \leq C.
\end{equation}
Since for that sequence we have $\frac{1}{\sqrt{E^h}} \leq C \frac{h^2}{E^h}$ we conclude that $$\| \frac{1}{\sqrt{E^h}}\sym ((\nabla \yb^h) \circ \Thetab^h\circ P^h-\Ibb) \|_{L^2(\Omega)} \leq C.$$ In the same way as before we conclude that $\| \frac{1}{\sqrt{E^h}} \sym \nabla' \ub^h \|_{L^2(\omega)}$ is bounded from which, again using Lemma \ref{Kornova}, (\ref{nacinzac}) and  (\ref{zazaKorna}), we conclude convergence (\ref{konvergencija2}). It remains to conclude $\Abb_{12}=0$. But this is easy, from (\ref{konvergencija2}) and (\ref{prilagodbba}), using the fact that $\min\{\frac{1}{\sqrt{E^h}},\frac{h^2}{E^h} \}\frac{\sqrt{E^h}}{h} \to \infty$ as $h \to 0$.
\end{prooof}
\begin{lemma} \label{naknadnododano}
Under the assumptions of Lemma \ref{osnovna} we can also conclude that
\begin{eqnarray*}
 &&\min\Big( \frac{h^2}{E^h}, \frac{1}{\sqrt{E^h}} \Big) \Big[(\widetilde{\yb}_{\alpha} \circ \Thetab^h \circ P^h-\Thetab^h_\alpha \circ P^h)+\sqrt{E^h} x_3\partial_{\alpha} v \Big] \rightharpoonup \ub_{\alpha}, \\
&&  \frac{h}{\sqrt{E^h}} (\widetilde{\yb}_3 \circ \Thetab^h \circ P^h-\Thetab_3^h \circ P^h) \to v,
\end{eqnarray*}
all in $W^{1,2}(\Omega)$.
\end{lemma}
\begin{prooof}
By the Korn's inequality and the relations (\ref{prilagodbba}), (\ref{nakkkkk}) it is enough to prove that
\begin{eqnarray} \nonumber \label{L2u}
o_{\alpha}^h &:=&\min\Big( \frac{h^2}{E^h}, \frac{1}{\sqrt{E^h}} \Big) \Big[(\widetilde{\yb}_{\alpha} \circ \Thetab^h \circ P^h-\Thetab^h_\alpha \circ P^h-\Ub^h_{\alpha})+\sqrt{E^h} x_3\partial_{\alpha} v \Big] \to 0,  \\ \label{oooostatku}
o_3^h &:=& \label{L2v} \frac{h}{\sqrt{E^h}} (\widetilde{\yb}_3 \circ \Thetab^h \circ P^h-\Thetab_3^h \circ P^h-V^h) \to 0,
\end{eqnarray}
both in $L^2(\Omega)$.
We shall just consider the situation $h^{-4}E^h \to \infty$ or $h^{-4}E^h \to 1$.
Using  Poincare  the inequality on each segment $\{x'\} \times [-\frac{1}{2},\frac{1}{2}]$, the Cauchy-Schwartz inequality and (\ref{prilagodbba}) we conclude
\begin{eqnarray} \nonumber
 \|o_{\alpha} ^h \|_{L^2(\Omega)} &\leq& C \frac{h^2}{E^h}\int_{\omega}  \|\partial_3 (\widetilde{\yb}_\alpha^h \circ \Thetab^h \circ P^h)(x',\cdot)-\partial_3(\Thetab^h_3 \circ P^h)(x',\cdot)\\ \nonumber & &\hspace{10ex}-\sqrt{E^h} \Abb_{\alpha 3}\|_{L^2(\{x'\} \times [-\frac{1}{2},\frac{1}{2}])}dx' \\ \nonumber &\leq & C \frac{h^2}{\sqrt{E^h}}\left\| \frac{h}{\sqrt{E^h}} \Big[\frac{1}{h} \partial_3 (\widetilde{\yb}_3^h \circ \Thetab^h \circ P^h)-\frac{1}{h} \partial_3(\Thetab^h_3 \circ P^h)\Big] -\Abb_{\alpha 3} \right\|_{L^2(\Omega)}\\ & & \to 0\\ \nonumber
 \|o_3^h \|_{L^2(\Omega)} &\leq& C \frac{h}{\sqrt{E^h}}\int_{\omega} \Big\|\partial_3 (\widetilde{\yb}_3^h \circ \Thetab^h \circ P^h)(x',\cdot)-\\ \nonumber & & \hspace{15ex} -\partial_3(\Thetab^h_3 \circ P^h)(x',\cdot)\Big\|_{L^2(\{x'\} \times [-\frac{1}{2},\frac{1}{2}])}dx' \\  &\leq & Ch \frac{h}{\sqrt{E^h}}\left\| \frac{1}{h} \partial_3 (\widetilde{\yb}_3^h \circ \Thetab^h \circ P^h)-\frac{1}{h} \partial_3(\Thetab^h_3 \circ P^h) \right\|_{L^2(\Omega)} \to 0.
\end{eqnarray}
\end{prooof}
\begin{lemma} \label{identifikacija}
Consider $\yb^h:\hat{\Omega}^h \to \ZR^3$, $\Rbb^h:\omega \to \SO(3)$ and $E^h>0$ and set
$$\ub^h:=\min \Big( \frac{h^2}{E^h},\frac{1}{\sqrt{E^h}} \Big) \int_{-1/2}^{1/2} (\left( \begin{array}{c} \yb^h_1 \circ \Thetab^h \circ P^h \\  \yb^h_2 \circ \Thetab^h \circ P^h \end{array} \right)(.,x_3)-x') dx_3,$$  $$
v^h:= \frac{h}{\sqrt{E^h}}\int_{-1/2}^{1/2} ((\yb_3^h \circ \Thetab^h\circ P^h)(.,x_3)-f^h \theta(.))  dx_3.  $$
Suppose that we have a subsequence of  $\yb^h$ such that
\begin{equation} \label{uvjetnaniz}  \lim_{h \to 0} h^{-2}E^h=0, \end{equation}
\begin{equation} \label{cudno1}
\| (\nabla \yb^h) \circ \Thetab^h\circ P^h -\Rbb^h \|_{L^2(\Omega)} \leq C\sqrt{E^h},
\end{equation}
\begin{equation} \label{cudno2}
\ub^h \rightharpoonup \ub \quad \textrm{in } W^{1,2}(\omega;\ZR^2), \quad v^h \to v \quad \textrm{in } W^{1,2}(\omega), \quad v \in W^{2,2}(\omega).
\end{equation}
Then
\begin{equation} \label{cudnoo}
\frac{h}{\sqrt{E^h}} (\Rbb^h-\Ibb) \to \Abb=\eb_3 \otimes \nabla' v-\nabla ' v \otimes \eb_3 \quad \textrm{in } L^2(\omega;\ZR^{3 \times 3}),
\end{equation}
and
\begin{equation}
\Gbb^h:= \frac{(\Rbb^h)^T ((\nabla \yb^h) \circ \Thetab^h\circ P^h)-\Ibb}{\sqrt{E^h}} \rightharpoonup \Gbb \quad \textrm{in } L^2(\Omega;\ZR^{3 \times 3}).
\end{equation}
and the $2 \times 2$ sub-matrix $\Gbb''$ given by $\Gbb''_{\alpha \beta}=\Gbb_{\alpha \beta}$ for $1 \leq \alpha,\beta \leq 2$ satisfies
\begin{equation} \label{zakor}
\Gbb''(x',x_3)=\Gbb_0(x')+x_3 \Gbb_1(x'),
\end{equation}
where
\begin{equation} \label{defg1}
\Gbb_1=-(\nabla')^2 v.
\end{equation}
Moreover
\begin{equation} \label{slucajmh4}
\nabla' \ub+(\nabla' \ub)^T+\nabla' v \otimes \nabla' v+\nabla' v \otimes \nabla' \theta+\nabla' \theta \otimes \nabla' v=0, \quad \textrm{if } h^{-4} E^h \to \infty,
\end{equation}
\begin{equation} \label{slucajh4}
\sym \Gbb_0=\frac{1}{2}(\nabla' \ub+(\nabla' \ub)^T+\nabla'v \otimes \nabla' v+\nabla' v \otimes \nabla' \theta+\nabla' \theta \otimes \nabla' v),  \quad \textrm{if } h^{-4}E^h \to 1,
\end{equation}
\begin{equation} \label{slucajvh4}
\sym \Gbb_0=\frac{1}{2}(\nabla' \ub+(\nabla' \ub)^T+\nabla' v \otimes \nabla' \theta+\nabla' \theta \otimes \nabla' v), \quad \textrm{if } h^{-4} E^h \to 0.
\end{equation}
\end{lemma}
\begin{prooof} We follow the proof of Lemma 15 in \cite{Muller3}.
We first assume (\ref{cudnoo}) (for the special sequence coming from Lemma \ref{osnovna} we know this anyhow) and establish the main assertion, namely the representation formula for $G$. Using the identity $2 \sym(\Qbb-\Ibb)=-(\Qbb-\Ibb)^T (\Qbb-\Ibb)$ which holds for all $\Qbb \in \SO(3)$ we immediately deduce from (\ref{cudnoo}) that
\begin{eqnarray} \nonumber
\frac{h^2}{E^h} \sym ( \Rbb^h-\Ibb) &\to& \frac{\Abb^2}{2}\\ \label{simetricnidio}   &=& -\frac{1}{2} \Big(\nabla'v \otimes \nabla'v+|\nabla ' v|^2  \eb_3 \otimes \eb_3 \Big) \ \textrm{in } L^1(\omega).
\end{eqnarray}
By the assumption $\Gbb^h$ is bounded in $L^2$, thus a subsequence converges weakly.
To show that the limit matrix $\Gbb''$ is affine in $x_3$ we consider the difference quotients
\begin{equation}
\Hbb^h(x',x_3)=s^{-1}[\Gbb^h(x',x_3+s)-\Gbb^h(x',x_3)].
\end{equation}
By multiplying the definition of $\Gbb^h$ with $\Rbb^h$ and using (\ref{kljucnozadokaz2}) we obtain for $\alpha,\beta \in \{1,2\}$
\begin{eqnarray} \nonumber
(\Rbb^h \Hbb^h)_{\alpha \beta} &=& \frac{1}{s \sqrt{E^h}} [\nabla_h (\yb^h \circ  \Thetab^h \circ P^h)(x',x_3+s) -\\ \nonumber & & \hspace{10ex}  -\nabla_h (\yb^h \circ  \Thetab^h \circ P^h)(x',x_3)]_{\alpha \beta} \\  \nonumber & &
+\frac{1}{s \sqrt{E^h}}\Big[((\nabla \yb^h) \circ \Thetab^h \circ P^h)(x',x_3+s)\cdot\\ \nonumber & & \hspace{10ex} \cdot \left(\Ibb-((\nabla \Thetab^h) \circ P^h)(x',x_3+s)\right)\Big]_{\alpha \beta} \\ \nonumber
& &-\frac{1}{s \sqrt{E^h}}\Big[((\nabla \yb^h) \circ \Thetab^h \circ P^h)(x',x_3) \cdot \\ & & \nonumber \hspace{10ex} \cdot \left(\Ibb-((\nabla \Thetab^h) \circ P^h)(x',x_3)\right)\Big]_{\alpha \beta} \\ \nonumber
&=&  \frac{1}{s \sqrt{E^h}} \Big[(\nabla_h (\yb^h \circ  \Thetab^h \circ P^h-\Thetab^h \circ P^h)(x',x_3+s)-\\ \nonumber & & \hspace{10ex}  -\nabla_h (\yb^h \circ  \Thetab^h \circ P^h-\Thetab^h \circ P^h)(x',x_3)\Big]_{\alpha \beta} \\ \nonumber
& &+\frac{1}{s \sqrt{E^h}}\Big[((\nabla \yb^h) \circ \Thetab^h \circ P^h)(x',x_3+s)-\\ \nonumber & & \hspace{10ex} -\Rbb^h(x'))\left(\Ibb-((\nabla \Thetab^h) \circ P^h)(x',x_3+s)\right)\Big]_{\alpha \beta} \\ \nonumber
& &-\frac{1}{s \sqrt{E^h}}\Big[((\nabla \yb^h) \circ \Thetab^h \circ P^h)(x',x_3)\\ \nonumber& &\hspace{10ex} -\Rbb^h(x'))\left(\Ibb-((\nabla \Thetab^h) \circ P^h)(x',x_3)\right)\Big]_{\alpha \beta}
 \\ \nonumber
& &-\frac{1}{s \sqrt{E^h}}\Big[(\Rbb^h (x')-\Ibb) ((\nabla \Thetab^h) \circ P^h)(x',x_3+s)-\\ & & \hspace{10ex} \label{zadnjarel}-((\nabla \Thetab^h) \circ P^h)(x',x_3)\Big]_{\alpha \beta}.
\end{eqnarray}
 For the first term in (\ref{zadnjarel}) we obtain
\begin{eqnarray}\nonumber
& &\frac{1}{s \sqrt{E^h}} [(\nabla_h (\yb^h \circ  \Thetab^h \circ P^h-\Thetab^h \circ P^h)(x',x_3+s)\\ \nonumber & & \hspace{10ex}-\nabla_h (\yb^h \circ  \Thetab^h \circ P^h-\Thetab^h \circ P^h)(x',x_3))]_{\alpha \beta}\\& & \hspace{7ex} =\frac{h}{\sqrt{E^h}} \partial_{\beta} \Big( \frac{1}{s} \int_0^s \frac{1}{h} \partial_3 (\yb^h \circ  \Thetab^h \circ P^h-\Thetab^h \circ P^h)_{\alpha} \Big).
\end{eqnarray}
From (\ref{cudno1}) and (\ref{cudnoo}) we conclude
\begin{equation}
 \frac{h}{\sqrt{E^h}}((\nabla \yb^h) \circ \Thetab^h \circ P^h-\Ibb) \to \Abb \ \textrm{in} \ L^2(\Omega;\ZR^{3 \times 3}).
\end{equation}
From that, using that $\|\nabla \Thetab^h -\Ibb\|_{L^\infty(\Omega^h)} \to 0$ and (\ref{kljucnozadokaz2}), we easily conclude that
\begin{equation} \label{velikpom}
\frac{h}{\sqrt{E^h}} \nabla_h (\widetilde{\yb}^h \circ \Thetab^h \circ P^h-\Thetab^h \circ P^h) \to \Abb \quad \textrm{in } L^2(\Omega;\ZR^{3 \times 3}).
\end{equation}
Using (\ref{velikpom}) we conclude that the first term in (\ref{zadnjarel}) converges weakly in $W^{-1,2}(\omega \times (-1,1-s))$ to $\Abb_{\alpha 3,\beta} (x')=-v_{,\alpha \beta} (x')$.
By using (\ref{glupaocjena}), (\ref{kristina}), (\ref{cudno1}) and (\ref{cudnoo}) we conclude that the all the rest converges to $0$ in $L^2(\omega \times (-1,1-s))$.
Since we have $\Rbb^h \to \Ibb$ boundedly a.e. and $\Hbb^h \rightharpoonup \Hbb$ in $L^2$  we thus obtain  $H_{\alpha \beta}(x',x_3)=-v_{,\alpha \beta}(x')$. From that we conclude that $\Gbb''$ is affine in $x_3$ and that $\Gbb_1$ has the form given in lemma. In order to prove formula for $\Gbb_0$ it suffices to study
$$ \Gbb_0^h (x')=\int_{-\frac{1}{2}}^{\frac{1}{2}} \Gbb^h(x',x_3) dx_3. $$
We have for $\alpha, \beta \in \{1,2\}$
\begin{eqnarray} \nonumber
(\Gbb^h)_{\alpha \beta}(x',x_3) &=&  \frac{((\nabla' \yb^h)\circ \Thetab^h \circ P^h-\Ibb)_{\alpha \beta}}{\sqrt{E^h}}-\frac{(\Rbb^h-\Ibb)_{\alpha \beta}}{\sqrt{E^h}}\\  \label{idensim}& &+\Big[(\Rbb^h-\Ibb)^T \frac{(\nabla \yb^h) \circ \Thetab^h \circ P^h-\Rbb^h}{\sqrt{E^h}} \Big]_{\alpha \beta}.
\end{eqnarray}
In the case $h^{-4}E^h \to 1$ using (\ref{nablaje}) , the identity (\ref{zaKorna}), (\ref{velikpom}) and the convergence of $\ubb^h$ we have
\begin{eqnarray} \label{dodatak1}
\sym  [\int_{-1/2}^{1/2} \frac{((\nabla' \yb^h)\circ \Thetab^h \circ P^h-\Ibb)}{\sqrt{E^h}}dx_3]_{\alpha \beta} \rightharpoonup
 \sym [\nabla ' \ub+ \Abb\Cbb]_{\alpha \beta},
\end{eqnarray}
in $L^2(\omega)$.
Using the convergence (\ref{simetricnidio}) we easily conclude
\begin{equation} \label{skorogot}
(\sym \Gbb_0^h)_{\alpha \beta} \rightharpoonup \Big[ \sym \nabla ' \ub+\sym (\Abb\Cbb)-\frac{\Abb^2}{2}\Big]_{\alpha \beta} \ \textrm{in } L^2(\omega).
\end{equation}
From (\ref{skorogot}) we easily deduce (\ref{slucajh4}).
To derive (\ref{slucajmh4}) let us observe that  in the case $h^{-4}E^h \to \infty$ we have $f^h=h^{-1}\sqrt{E^h}$. Multiplying (\ref{idensim}) by  $h^2/\sqrt{E^h}$ ($\to 0$) we
conclude
\begin{equation}\frac{h^2}{E^h}  [\sym\int_{-1/2}^{1/2} ((\nabla' \yb^h)\circ \Thetab^h \circ P^h-\Ibb) dx_3]_{\alpha \beta}- \sym[\frac{h^2}{E^h} (\Rbb^h-\Ibb)]_{\alpha \beta} \to 0 \ \textrm{in } L^1(\omega).
\end{equation}
Using again (\ref{nablaje}) , the identity (\ref{zaKorna}) and the convergence of $\ubb^h$ we have
\begin{equation}
\frac{h^2}{E^h}  [\sym\int_{-1/2}^{1/2} ((\nabla' \yb^h)\circ \Thetab^h \circ P^h-\Ibb) dx_3]_{\alpha \beta} \rightharpoonup \sym [\nabla ' \ub+ \Abb\Cbb]_{\alpha \beta}  \ \textrm{in } L^2(\omega).
\end{equation}
Using the convergence (\ref{simetricnidio})  we conclude (\ref{slucajmh4}).
(\ref{slucajvh4}) can be concluded in the similar way from (\ref{idensim}).
It remains to prove (\ref{cudnoo}). Since $\Rbb^h$ is independent of $x_3$ we have for $i,j  \in \{1,2,3\}$
\begin{equation} \label{zaborav1}
(\Rbb^h - \Ibb)_{i j} =\int_{-1/2}^{1/2} (\Rbb^h-(\nabla \yb^h) \circ \Thetab^h\circ P^h)_{i j}dx_3+\int_{-1/2}^{1/2} ( ( \nabla \yb^h) \circ \Thetab^h\circ P^h-\Ibb)_{i j}dx_3 \\ \label{zaborav2}.
\end{equation}
By using (\ref{glupaocjena}), the relation (\ref{zaKorna}) and the convergence (\ref{cudno2}) we conclude for $\alpha,\beta \in \{1,2\}$
\begin{equation} \label{ocjenazaunutra}
\| \int_{-1/2}^{1/2} ( (\nabla \yb^h) \circ \Thetab^h\circ P^h-\Ibb)_{\alpha \beta}dx_3\|_{L^2(\omega)} \leq C(h^{-1} \sqrt{E^h} f^h+ \max\{\sqrt{E^h},\frac{E^h}{h^2}\}).
\end{equation}
From (\ref{zaborav1}) and (\ref{ocjenazaunutra}), by using  (\ref{cudno1}), we conclude
\begin{equation} \label{strpljenje}
\frac{h}{\sqrt{E^h}} (\Rbb^h-\Ibb)_{\alpha \beta} \to 0 \ \textrm{in } L^2(\omega).
\end{equation}
In the same way, by using (\ref{glupaocjena}), the relation (\ref{zaKorna}),  (\ref{cudno1}) and the convergence (\ref{cudno2}), we conclude for $\beta \in \{1,2\}$
\begin{eqnarray*} \nonumber
& &\| \int_{-1/2}^{1/2} ( (\nabla \yb^h) \circ \Thetab^h\circ P^h-\Ibb)_{3 \beta}dx_3- (h^{-1} \sqrt{E^h}) \partial_{\beta} v\|_{L^2(\omega)}\\ \label{ocjenadonjih} & & \hspace{15ex} \leq C(h^{-1} \sqrt{E^h} f^h+ \max\{\sqrt{E^h},\frac{E^h}{h^2}\}).
\end{eqnarray*}
In the same way as for (\ref{strpljenje})  we conclude
\begin{equation} \label{strpljenje2}
\frac{h}{\sqrt{E^h}} (\Rbb^h-\Ibb)_{3 \beta} \to \partial_{\beta} v \ \textrm{in } L^2(\omega).
\end{equation}
Using the fact that $\Rbb^h$ takes the values  in $\SO(3)$ we deduce that $\|\Rbb^h_{\beta 3}\|_{L^2(\omega)} \leq C \sqrt{E^h}/h$.
To get control on $\Rbb_{33}^h$ we use the fact that for $\Qbb \in \SO(3)$ we have
$$|1-\Qbb_{33}|=| \det \Qbb-\Qbb_{33}| \leq C \sum_{\alpha, \beta=1}^2 |(\Qbb-\Ibb)_{\alpha \beta}|+ C(|\Qbb_{13} \Qbb_{31}|+|\Qbb_{23}\Qbb_{32}|). $$
From this, using the generalized convergence theorem (with $L^2$ convergent majorant rather than constant majorant), we easily deduce that $(h/\sqrt{E^h} (\Rbb^h_{33}-1) \to 0$ in $L^2(\omega)$. To control $\Rbb^h_{13}$ we use the fact that the first and third row of $\Rbb^h$ are orthogonal. This yields
$$ | \Rbb^h_{13} +\Rbb^h_{31} | \leq C(|\Rbb_{11}^h-1|+|\Rbb^h_{33}-1|+|\Rbb^h_{12}|), $$
and together with (\ref{strpljenje}), (\ref{strpljenje2}) and the convergence of $\Rbb^h_{33}$ this gives the desired convergence for $\Rbb_{13}^h$. The same argument applies to $\Rbb_{23}^h$ and this finishes the proof.
\end{prooof}

By $I^{MvK}$ and $I^{LMvK}$ we denote the functionals
\begin{eqnarray} \nonumber
I^{MvK}(\ub,v)&=&\int_{\omega} \Big( \frac{1}{2} Q_2(\frac{1}{2}(\nabla' \ub+(\nabla' \ub)^T+\nabla' v \otimes \nabla' v+\nabla' v \otimes \nabla' \theta\\ & & \hspace{10ex} +\nabla' \theta \otimes \nabla' v)  \label{defMvK}
+ \frac{1}{24}Q_2((\nabla ')^2 v) \Big) dx',
\\ \nonumber
I^{LMvK}(\ub,v)&=&\int_{\omega} \Big( \frac{1}{2} Q_2(\frac{1}{2}(\nabla' \ub+(\nabla' \ub)^T+\nabla' v \otimes \nabla' \theta+\nabla' \theta \otimes \nabla' v) \\  \label{defLMvK}
&+& \frac{1}{24}Q_2((\nabla ')^2 v) \Big) dx',
\end{eqnarray}
defined on the space $W^{1,2}(\omega; \ZR^2) \times W^{2,2}(\omega)$.
\begin{remark}
The term $\nabla '\ub+(\nabla' \ub)^T+\nabla' v \otimes \nabla' v+\nabla' v \otimes \nabla' \theta+\nabla' \theta \otimes \nabla' v$ measures the change of the metric tensor (which is of the second order), while the term $(\nabla ')^2 v$ measures the curvature tensor change (which is of the first order) of the deformation $\varphi(x') = \left( \begin{array}{c} x'+(f^h)^2 \ub (x') \\ f^h (v+\theta)(x') \end{array} \right)$ with respect to the deformation $\varphi_0(x') = \left( \begin{array}{c} x' \\ f^h \theta(x') \end{array} \right)$. To see this let us calculate
\begin{eqnarray*}
(\nabla' \varphi)^T \nabla' \varphi &=& \Big[ \left( \begin{array}{c}\Ibb \\ 0  \end{array} \right)+ \left( \begin{array}{c} (f^h)^2 \nabla' \ub  \\ f^h  \nabla' (v+\theta)(x') \end{array} \right) \Big]^T \Big[ \left( \begin{array}{c}\Ibb \\ 0  \end{array} \right)+ \left( \begin{array}{c} (f^h)^2 \nabla' \ub  \\ f^h  \nabla' (v+\theta)(x') \end{array} \right) \Big] \\
&=& \Ibb+ (f^h)^2 (\nabla '\ub+(\nabla' \ub)^T+\nabla' (v+\theta) \otimes \nabla' (v+\theta))+ O((f^h)^4)
\\
&=& (\nabla' \varphi_0)^T \nabla' \varphi_0 +(f^h)^2 (\nabla '\ub+(\nabla' \ub)^T+\nabla' v \otimes \nabla' v+\nabla' v \otimes \nabla' \theta\\ & & \hspace{25ex}+\nabla' \theta \otimes \nabla' v) + O((f^h)^4).
\end{eqnarray*}
The same explanation goes in the situations when $h^{-4} E^h \to 0$ when the term $\nabla' v \otimes \nabla' v$ disappears.
To analyze the change of the curvature we need to calculate the normals
\begin{eqnarray*}
\nb_{\varphi} &=& \eb_3-(f^h) \partial_1 (v+\theta) \eb_1-(f^h) \partial_2  (v+\theta) \eb_2+O((f^h)^2), \\
\nb_{\varphi_0} &=& \eb_3-(f^h) \partial_1 \theta \eb_1-(f^h) \partial_2  \theta \eb_2+O((f^h)^2).
\end{eqnarray*}
Now we have for arbitrary $\alpha, \beta \in \{1,2\}$
\begin{eqnarray*}
\partial_{\alpha \beta} \varphi \cdot \nb_{\varphi} &=& f^h (\partial_{\alpha \beta} (v+\theta))+O((f^h)^2) \\
&=&  \partial_{\alpha \beta} \varphi_0 \cdot \nb_{\varphi_0}+f^h \partial_{\alpha \beta} v+O((f^h)^2)
\end{eqnarray*}
\end{remark}
\begin{corollary} \label{zajakuisto}
Let $E^h, \yb^h,\Rbb^h,\ub^h, v^h, \Gbb^h, \Gbb, \Gbb'', \Gbb_0, \Gbb_1$ be as in  Lemma \ref{identifikacija}. Then we have the following semi-continuity results.
\begin{enumerate}[i)]
\item If $\lim_{h \to 0} h^{-4} E^h=\infty$ then
\begin{equation}
\liminf_{h \to 0} \frac{1}{E^h} I^h (\yb^h) \geq \int_{\omega} \frac{1}{24} Q_2 ((\nabla')^2 v) dx'.
\end{equation}
\item If $\lim_{h \to 0} h^{-4}E^h=1$ then
\begin{eqnarray}
\nonumber \liminf_{h \to 0} \frac{1}{E^h} I^h(\yb^h) &\geq& I^{MvK}(\ub,v).
\end{eqnarray}
\item If  $\lim_{h \to 0} h^{-4}E^h=0$ then
\begin{eqnarray}
\nonumber \liminf_{h \to 0} \frac{1}{E^h} I^h(\yb^h) &\geq& I^{LMvK}(\ub,v).
\end{eqnarray}
\end{enumerate}
\end{corollary}
\begin{prooof}
 We shall use the truncation, the Taylor expansion and the weak semi-continuity argument as in the proof of Corollary $16$ in \cite{Muller3}.
Let $m:[0,\infty) \to [0,\infty)$ denote a modulus of continuity of $D^2W$ near the identity and consider the good set $\Omega_h :=\{x \in \Omega: |\Gbb^h (x)| < h^{-1} \}$. Its characteristic function $\chi_h $ is bounded and satisfies $\chi_h \to 1$ a.e. in $\Omega$. Thus we have $\chi_h \Gbb^h \rightharpoonup \Gbb$ in $L^2 (\Omega)$. By Taylor expansion
\begin{equation}
\frac{1}{E^h} \chi_h W(\Ibb+\sqrt{E^h} \Gbb^h) \geq \frac{1}{2} Q_3(\chi_h \Gbb^h)-m(h^{-1} \sqrt{E^h})|\Gbb^h|^2.
\end{equation}
Using (\ref{uvjetnaniz}) and the boundedness of  sequence $\Gbb^h$ in $L^2(\Omega;\ZR^{3 \times 3})$ we conclude
\begin{eqnarray}
\nonumber & & \liminf_{h \to 0} \frac{1}{E^h} I^h(\yb^h) \\
\nonumber &=& \liminf_{h \to 0} \frac{1}{E^h} \int_{\Omega} W((\Rbb^h)^T ((\nabla \yb^h) \circ \Thetab^h\circ P^h)) dx \\
\nonumber & \geq & \Big[ \frac{1}{2} \int_\Omega Q_3(\chi_h \Gbb^h) dx + \frac{1}{E^h} \int_\Omega (1-\chi_h) W((\nabla \yb^h) \circ \Thetab^h\circ P^h)) dx \Big] \\
& \geq & \frac{1}{2} \int_{\Omega} Q_3 (\Gbb) dx \geq \frac{1}{2} \int_{\Omega} Q_2 (\Gbb'') dx.
\end{eqnarray}
Here we have used the fact that $Q_3$ is a positive semi-definite quadratic form and therefore the functional $v \mapsto \int_\Omega Q_3 (v)$ is weakly lower semi-continuous in $L^2$. Now by (\ref{zakor}) we have
\begin{equation}
\int_{-1/2}^{1/2} Q_2 (\Gbb'')(x',x_3) dx_3 = Q_2 (\Gbb_0(x'))+ \frac{1}{12} Q_2 (\Gbb_1 (x')).
\end{equation}
 Together with (\ref{defg1}), (\ref{slucajh4}), (\ref{slucajvh4}) this implies the claim of the corollary.
\end{prooof}

\subsection{Upper bound}
\begin{theorem}[optimality of lower bound] \label{upperbound}
If $h^{-4}E^h \to 1$ and if  $v \in W^{2,2}(\omega)$,  $\ub \in W^{1,2} (\omega; \ZR^2)$ then there exists a sequence
$\hat{\yb}^h \in W^{1,2}(\hat{\Omega}^h;\ZR^3)$ such that
\begin{equation} \label{devide}
(\nabla \hat{\yb}^h) \circ \Thetab^h \circ P^h \to\Ibb \ \textrm{in } L^2(\Omega; \ZR^{3 \times 3}),
 \end{equation}
 and for $\Ub^h$, $V^h$ defined by (\ref{definicijauv}) (where $\widetilde{\yb}$ should be replaced by $\hat{\yb}$)  convergence (\ref{konvergencija})-(\ref{konvergencija2}) are valid
and
\begin{equation}
\lim_{h \to 0} \frac{1}{E^h} I^h (\hat{\yb}^h)= I^{MvK}(\ub,v).
\end{equation}
If $h^{-4}E^h \to 0$ and if  $v \in W^{2,2}(\omega)$,  $\ub \in W^{1,2} (\omega; \ZR^2)$ then there exists
$\hat{\yb}^h$ such that the convergence (\ref{devide}),  (\ref{konvergencija})-(\ref{konvergencija2}) hold and
\begin{equation} \label{konvlmvk}
\lim_{h \to 0} \frac{1}{E^h} I^h (\hat{\yb}^h)= I^{LMvK}(\ub,v).
\end{equation}
\end{theorem}
\begin{prooof}
Let us first analyze the situation when $h^{-4}E^h \to 1$ and assume that $\ub,v$ are smooth. Then we define
\begin{eqnarray} \nonumber
\hat{\yb}^h (\Thetab^h(x',x_3^h)) &=& 
\Thetab^h (x',x_3)+\left(\begin{array}{c} \frac{E^h}{h^2} \ub (x') \\ \frac{\sqrt{E^h}}{h} v(x') \end{array} \right)-\frac{\sqrt{E^h}}{h} x_3^h \left(\begin{array}{c} \partial_1 v(x') \\ \partial_2 v(x') \\ 0 \end{array} \right)+\\ \label{defhaty} & &+ \frac{E^h}{h^2}x_3^h \db^0(x')+\frac{1}{2} \frac{\sqrt{E^h}}{h} (x_3^h)^2 \db^1 (x'),
\end{eqnarray}
where $\db_0,\db_1: \omega \to \ZR^3$ are going to be chosen later.
The convergence (\ref{devide}) as well the convergences  (\ref{konvergencija})-(\ref{konvergencija2}) can easily seen to be valid for this sequence. We also have
\begin{eqnarray} \nonumber
\nabla \hat{\yb}^h \nabla \Thetab^h &=&\nabla \Thetab^h + \frac{\sqrt{E^h}}{h}
  \left( \begin{array}{c|c} 0 & - (\nabla ' v)^T \\ \hline  \nabla ' v & 0
  \end{array} \right)+
  \frac{E^h}{h^2}\left( \begin{array}{c|c} \nabla' \ub & 0 \\ \hline  0 & 0
  \end{array} \right) \\ \nonumber  & &- \frac{\sqrt{E^h}}{h} x_3^h \left( \begin{array}{c|c}  (\nabla')^2 v & 0 \\ \hline  0 & 0
  \end{array} \right)
  + \frac{E^h}{h^2} \db^0 \otimes \eb_3 +  \frac{\sqrt{E^h}}{h} x_3 ^h \db^1 \otimes \eb_3 \\ & & \label{nablayt} +O(h^3).
\end{eqnarray}
From (\ref{nablayt}) by using (\ref{inverz}) we conclude
\begin{eqnarray} \nonumber
\nabla \hat{\yb}^h &=& \Ibb+\frac{\sqrt{E^h}}{h}
  \left( \begin{array}{c|c} 0 & - (\nabla ' v)^T \\ \hline  \nabla ' v & 0
  \end{array} \right)+
  \frac{E^h}{h^2}\left( \begin{array}{c|c} \nabla' \ub+\nabla' v \otimes \nabla' \theta & 0 \\ \hline  0 & 0
  \end{array} \right) \\ \nonumber & &- \frac{\sqrt{E^h}}{h} x_3^h \left( \begin{array}{c|c}  (\nabla')^2 v & 0 \\ \hline  0 & 0
  \end{array} \right)
  + \frac{E^h}{h^2} \db^0 \otimes \eb_3 +  \frac{\sqrt{E^h}}{h} x_3 ^h \db^1 \otimes \eb_3 \\ & & \label{nablayyt} +O(h^3).
\end{eqnarray}
Using the identities $(\Ibb+\Abb)^T (\Ibb+\Abb)=\Ibb+2 \sym \Abb +\Abb^T \Abb$ and $(\eb_3 \otimes \ab'-\ab' \otimes \eb_3)^T (\eb_3 \otimes \ab'-\ab' \otimes \eb_3)=\ab' \otimes \ab' +|\ab'|^2 \eb_3 \otimes \eb_3$ for $\ab' \in \ZR^2$ we obtain
\begin{eqnarray} \nonumber
(\nabla \hat{\yb}^h)^T (\nabla \hat{\yb}^h) &=& \Ibb+ \frac{E^h}{h^2}[ 2\sym (\nabla ' \ub+\nabla' v \otimes \nabla'\theta)+ \nabla' v \otimes \nabla 'v+ \\ \nonumber& &\hspace{1ex} +|\nabla ' v|^2 \eb_3 \otimes \eb_3] -2\frac{\sqrt{E^h}}{h} x_3^h (\nabla')^2 v+
2 \frac{E^h}{h^2} \sym (\db^0 \otimes \eb_3)+\\ & &+2\frac{\sqrt{E^h}}{h} x_3 ^h \sym(\db^1 \otimes \eb_3) +O(h^3).
\end{eqnarray}
Taking the square root and using the frame indifference of $W$ and the Taylor expansion we get
\begin{equation}
(E^h)^{-1} W(\nabla  \hat{\yb}^h)=(E^h)^{-1} W([(\nabla  \hat{\yb}^h)^T \nabla  \hat{\yb}^h ]^{1/2}) \to \frac{1}{2} Q_3 (\Abb+x_3 \Bbbb),
\end{equation}
where
\begin{eqnarray*}
\Abb &=& \sym(\nabla ' \ub+ \nabla ' v \otimes \nabla ' \theta)+\frac{1}{2} \nabla ' v \otimes \nabla ' v + \frac{1}{2} | \nabla' v|^2 \eb_3 \otimes \eb_3 + \sym (\db^0 \otimes \eb_3), \\
\Bbbb &=& - (\nabla ' )^2 v + \sym (\db^1 \otimes \eb_3).
\end{eqnarray*}
For a symmetric $2 \times 2 $ matrix $\Abb''$ let $c=\mathcal{L} \Abb'' \in \ZR^3$ denote the (unique) vector which realizes the minimum in the definition of  $Q_2$. i.e.
$$ Q_2 (\Abb'')=Q_3 (\Abb''+\cb \otimes \eb_3+\eb_3 \otimes \cb). $$
Since $Q_3$ is positive definite on symmetric matrices, $\cb$ is uniquely determined and the map $\mathcal{L}$ is linear. We now take
\begin{eqnarray*}
\db^0 &=& -\frac{1}{2} | \nabla ' v|^2 \eb_3 + \mathcal{L} (\sym(\nabla ' \ub+ \nabla ' v \otimes \nabla ' \theta)+ \nabla ' v \otimes \nabla ' v), \\
\db^1 &=& -2 \mathcal{L} ((\nabla')^2 v).
\end{eqnarray*}
This finishes the proof of theorem in the situation $h^{-4}E^h \to 1$ and smooth $\ub, v$. For general $\ub,v$ it suffices to consider  smooth approximations $\ub^h , v^h, \db^{0,h}, \db^{1,h}$  of $\ub,v$ in $W^{1,2}(\omega)$ i.e. $W^{2,2}(\omega;\ZR^2)$. We first choose $\ub^h \in C^\infty(\omega;\ZR^2), v^h \in C^\infty (\omega)$ such that
$\| \ub^h-\ub\|_{W^{1,2}} < h$ and $\| v^h-v\|_{W^{2,2}} < h$ and $|I^{MvK}(\ub^h,v^h)-I^{MvK} (\ub, v)|<h$. For $\ub^h, v^h$ we choose $\hat{\yb}^h$ such that $$|\frac{1}{E^h} I^h(\hat{\yb}^h)-I^{MvK}(\ub^h,v^h)|< h.$$ Then we have the claim.

In the situation when $h^{-4}E^h \to 0$, for smooth $\ub,v$, we define
\begin{eqnarray} \nonumber
\hat{\yb}^h (\Thetab^h(x',x_3^h)) &=& 
\Thetab^h (x',x_3)+\left(\begin{array}{c} \sqrt{E^h} \ub (x') \\ \frac{\sqrt{E^h}}{h} v(x') \end{array} \right)-\frac{\sqrt{E^h}}{h} x_3^h \left(\begin{array}{c} \partial_1 v(x') \\ \partial_2 v(x') \\ 0 \end{array} \right)+\\& &+ \sqrt{E^h} x_3^h \db^0(x')+\frac{1}{2} \frac{\sqrt{E^h}}{h} (x_3^h)^2 \db^1 (x').
\end{eqnarray}
In the same way as in the situation $h^{-4}E^h \to 1$ we can define $\db^0$, $\db^1$ such that (\ref{konvlmvk}) is satisfied (the term $\nabla' v \otimes \nabla ' v$ disappears because it is of the order $\frac{E^h}{h^2}$ for which  $\frac{1}{\sqrt{E^h}} \frac{E^h}{h^2}\to 0$ is satisfied. For non-smooth $\ub, v$ the argument is the same as before.

\end{prooof}
\begin{remark}
The equality (\ref{slucajmh4}) can be written in the form
\begin{equation}
\nabla' \ub+(\nabla' \ub)^T=\nabla' \theta \otimes \nabla' \theta-\nabla' (v+\theta) \otimes \nabla' (v+\theta). \end{equation}
The left hand side is the symmetrized gradient of a $W^{1,2}$ function. It is known fact that, if $\omega$ is simply connected, a $L^2$ map $e: \omega \to \ZR^{2 \times 2}_{\sym}$ is the symmetrized gradient of a $W^{1,2}(\omega; \ZR^2)$ function $\ub$, i.e.
\begin{equation} 2e=(\nabla ' \ub)^T + \nabla ' \ub, \end{equation}
if and only if
\begin{equation}
\partial_{22} e_{11}+\partial_{11} e_{22}-2 \partial_{12} e_{12}=0,
\end{equation}
in the sense of distributions.
On the other hand  it is easily seen that this condition on the right hand side implies  (see Proposition 30 in \cite{Muller3})
\begin{equation} \label{uvjetjednakosti}
\det (\nabla')^2 (v+\theta)=\det (\nabla')^2 \theta.
\end{equation}
For smooth $v$ this means that the graphs of the functions $x' \to (v+\theta)(x')$, $x' \to \theta(x')$ have equal Gauss-Kronecker curvature at each point which is a necessary condition for the existence of an exact isometry between these two surfaces.
\end{remark}
\begin{remark} \label{napomenaoost}
Situation $h^{-2}E^h \to 0$ and $\lim_{h \to 0} h^{-4} E^h \to
\infty$ remains uncovered for the upper bound. There are
several reasons for that and they are similar to the one
observed in \cite{Muller3,Lewicka1,Lewicka2,Lewicka3}. In the
situation $h^{-3} E^h \to 0$ we need some additional regularity
results (e.g. $v \in W^{1,\infty}(\omega)$). This can be
concluded when the graph of $\theta$ is developable surface
i.e. when we have $\det (\nabla')^2 \theta=0$ (see
\cite{Muller3}). The situation $h^{-3} E^h >0$ is more
complicated. For given $v \in W^{2,2}(\omega)$ we would like to
construct the (exact or higher precision-see \cite{Lewicka3})
isometry (from the graph of the function $x' \to f^h
\theta(x')$ to the graph of the function  $x' \to f^h
(v+\theta)(x')$)  of the form
\begin{equation}
\bar{\yb}_{h}: \omega \to \ZR^3, \quad \bar{\yb}_{h}(x')=\left( \begin{array}{c} x'+(f^h)^2 \ub_h (x') \\ f^h (v+\theta)(x') \end{array} \right).
\end{equation}
The condition (\ref{uvjetjednakosti}) is a necessary condition for the existence of isometry, but not sufficient unless we have special situation $\det^2 (\nabla' \theta)=K$.
In the situation $K=0$ one can construct an exact isometry in the similar way as in \cite{Muller3} (see also \cite{Muller6}).
This can be done under additional (mild) hypothesis that there exists $\epsilon>0$ such that $h^{-(2+\epsilon)}E^h \to 0$.
The construction of an isometry would go by these steps:
\begin{enumerate}
\item it can be seen that there exists the
    $C^2(\bar{\omega})$ isometry between $\bar{\omega}$ and
    the graph of $f^h\theta(\bar{\omega})$ of the form
    $$i(x')=\left( \begin{array}{c} x'+(f^h)^2 \phi_h (x')
    \\ f^h \theta(x') \end{array} \right),$$ where $\phi_h
    \in C^2(\bar{\omega})$ (see the proof of Theorem 25 in
    \cite{Muller3}).
\item Using this isometry and the one between $\omega$ and
    the graph of $f^h (v+\theta) (\omega)$, one can easily
    construct an isometry between the graph of $f^h
    \theta(\omega)$ and the graph of $f^h(v+\theta)
    (\omega)$.
\end{enumerate}
In the general case, when $\theta(\omega)$ is not developable surface, in the situations $h^{-2}E^h \to 0$ and $h^{-3}E^h \to \infty$, stronger influence of the geometry of $\theta(\omega)$ on the model is expected (see \cite{Lewicka3}).
\end{remark}
\subsection{Convergence theorem}
Let $\fb_3^h \in L^2(\hat{\Omega}^h;\ZR)$ be given with the property
\begin{equation} \label{pretpostavkesile}
\int_{\hat{\Omega}^h} \fb_ 3^h=0, \ \quad \frac{1}{h\sqrt{E^h}}\fb_3^h \circ \Thetab^h \circ P^h \rightharpoonup \fb_3 \ \textrm{in } L^2(\Omega;\ZR).
\end{equation}
Let $m^h$ be the maximized action of force $\fb_3^h$ over all rotations of $\hat{\Omega}^h$,
\begin{equation} m^h= \max_{\Qbb \in \SO(3)}\int_{\hat{\Omega}^h} \fb_3^h(x) (\Qbb x)_3 dx, \end{equation}
and define
\begin{equation}
\mathcal{M}=\{ \bar{\Qbb} \in \SO(3); r(\bar{\Qbb})<+\infty \},
\end{equation}
to be the effective domain of the following relaxation functional $r:\SO(3) \to [0,+\infty]$
\begin{equation} \label{defodr}
r(\Qbb)= \min \Big\{ \liminf \frac{1}{hE^h} \Big( m^h-\int_{\hat{\Omega}^h} \fb^h_3(x) (\Qbb^h x)_3 dx\Big); \Qbb^h \in \SO(3), \Qbb^h \to \Qbb \Big\}.
\end{equation}
The set $\mathcal{M}$ identifies the candidates for large
rotation the body would perform to reduce its energy (see
\cite{Lewicka1,Lewicka2}).
\begin{remark}
The set $\mathcal{M}$ is introduced in \cite{Lewicka2}.  We can conclude (using the uniform convergence of the functions $\Qbb \to \frac{1}{h^2 \sqrt{E^h}}  \int_{\hat{\Omega}^h} \fb^h_3(x) (\Qbb x)_3 dx $) that $\frac{m^h} {h^2 \sqrt{E^h}}\to m$ where
\begin{eqnarray*} m &=& \max_{\Qbb \in \SO(3)}\int_{\omega} \Big( \int_{-1/2}^{1/2}\fb_3(x',x_3)dx_3 \Big)   (\Qbb \left( \begin{array}{c} x' \\ 0 \end{array} \right))_3 dx' \\ &=& \max_{q_1^2+q_2^2 \leq 1}\int_{\omega}\Big( \int_{-1/2}^{1/2}\fb_3(x',x_3)dx_3 \Big) (q_1 x_1+q_2 x_2) dx'.
\end{eqnarray*}
From that and the definition of $\mathcal{M}$ and $r$ it can be concluded that $\mathcal{M}\subset\mathcal{M}_0$
where
$$ \mathcal{M}_0 = \{ \bar{\Qbb} \in \SO(3): \int_{\omega}\Big( \int_{-1/2}^{1/2}\fb_3(x',x_3)dx_3 \Big) \ (\bar{\Qbb}  \left( \begin{array}{c} x' \\ 0 \end{array} \right))_3 dx'= m \}.$$
Since $r$ is lower semi-continuous it can be concluded that $\mathcal{M}$ is a nonempty closed subset of $\mathcal{M}_0$. Since we also have the change of geometry as $h \to 0$ it can not be concluded that $\mathcal{M} =\mathcal{M}_0$ even under the condition $  \frac{1}{h\sqrt{E^h}}\fb_3^h \circ \Thetab^h \circ P^h = \fb_3$ for every $h>0$. Let us mention that for every $\bar{\Qbb} \in \mathcal{M}_0 \supset \mathcal{M}$ we have the following equality
\begin{equation} \label{jednakostmomenata}
\int_{\omega} \Big(\int_{-1/2}^{1/2}\fb_3(x',x_3)dx_3\Big)  (\bar{\Qbb}\Fbb\left( \begin{array}{c} x' \\ 0 \end{array} \right))_3 dx'=0, \quad \forall \Fbb \in \so(3)
\end{equation}
The equality (\ref{jednakostmomenata}) is the consequence of
the fact that the differential vanishes at the extreme points
and the fact that $\so(3)$ is tangential to $\SO(3)$. The
equality (\ref{jednakostmomenata}) is the balance of momentum.
Since the "dead loads" $\fb_3$ are given on the reference
configuration, the shell adjusts its deformation to satisfy the
balance of momentum.
\end{remark}
To the total energy functional (defined on the space $W^{1,2}(\hat{\Omega}^h;\ZR^3)$) we add the constant and  redefine
\begin{equation} \label{defjh}
J^h (\yb^h)= I^h(\yb^h)- \frac{1}{h} \int_{\hat{\Omega}^h} \fb_3^h \yb_3^h +\frac{m^h}{h}.
\end{equation}
The following theorem is the main result.
\begin{theorem}[$\Gamma$-convergence] \label{najglavnijiteorem}
Assume $h^{-4}E^h \to 1$ or $h^{-4}E^h \to 0$. Let us suppose that $\fb_3^h \in L^2(\hat{\Omega}^h;\ZR)$ is given and satisfies (\ref{pretpostavkesile}). Then:
\begin{enumerate}[1.]
\item  There exists $C>0$ such that for every $h>0$ we have
\begin{equation} \label{nejinf} 0 \geq \inf\left\{ \frac{1}{E^h} J^h (\yb^h) ; \yb^h \in W^{1,2} ( \hat{\Omega}^h;\ZR^3) \right\} \geq -C.   \end{equation}
\item If $\yb^h \in W^{1,2}(\hat{\Omega}^h;\ZR^3)$ is a minimizing sequence of $\frac{1}{E^h} J^h$, that is
\begin{equation} \label{svojmin} \lim_{h \to 0} \Big( \frac{1}{E^h} J^h (\yb^h)- \inf \frac{1}{E^h} J^h \Big)=0, \end{equation}
then we have that there exists $\bar{\Rbb}^h \in \SO(3), \ \cb^h \in \ZR$ such that  the sequence $(\bar{\Rbb}^h, \yb^h)$ has its subsequence (also not relabeled) with the following property:
\begin{enumerate}[i)]
\item For the sequence $ \widetilde{\yb}^h:=(\bar{\Rbb}^h)^T \yb^h- \cb^h$ and $\Ub^h, V^h$ defined by (\ref{definicijauv}) the following is valid
    \begin{eqnarray*}
v^h &:=&\frac{h}{\sqrt{E^h}} V^h \to v  \quad \textrm{in } W^{1,2}(\omega),\  v \in W^{2,2}(\omega), \\
\ub^h &:=& \min\Big( \frac{h^2}{E^h}, \frac{1}{\sqrt{E^h}} \Big) \Ub^h \rightharpoonup \ub \quad \textrm{weakly in } W^{1,2}(\omega;\ZR^2).
\end{eqnarray*}
\end{enumerate}
Also, any accumulation point $\bar{\Rbb}$ of the sequence $\Rbb^h$ belongs to $\mathcal{M}$.
Moreover if $h^{-4}E^h \to 1$ then any accumulation point $(\ub,v,\bar{\Rbb})$ of the sequence $(\ub^h,v^h,\bar{\Rbb}^h)$ minimizes the functional
\begin{equation}
J^0(\ub,v,\bar{\Rbb})=I^{MvK}(\ub,v)- \bar{\Rbb}_{33} \int_{\omega} \Big(\int_{-1/2}^{1/2} \fb_3 dx_3 \Big)  v(x')dx'+r(\bar{\Rbb}),
\end{equation}
where $I^{MvK}$ is defined in (\ref{defMvK}) and $r$ is
defined in (\ref{defodr}). If $h^{-4}E^h \to 0$ then we
 have that any accumulation
point $(\ub,v,\bar{\Rbb})$ of the sequence
$(\ub^h,v^h,\bar{\Rbb}^h)$ minimizes the functional
\begin{equation}
J^0_L(\ub,v,\bar{\Rbb})=I^{LMvK}(\ub,v)- \bar{\Rbb}_{33} \int_{\omega} \Big(\int_{-1/2}^{1/2} \fb_3 dx_3 \Big)  v(x')dx'+r(\bar{\Rbb}),
\end{equation}
where $I^{LMvK}$ is defined in (\ref{defLMvK}).
\item The minimum of the functional $J^0$ i.e. $J^0_L$
    exists in the space $W^{1,2}(\omega;\ZR^2) \times
    W^{2,2}(\omega) \times \SO(3)$. If $\yb^h \in
    W^{1,2}(\hat{\Omega}^h;\ZR^3)$ is a minimizing sequence
    (not relabeled) of $\frac{1}{E^h} J^h$ then we have
    that
\begin{equation}
    \lim_{h \to 0} \frac{1}{E^h} J^h (\yb^h)= \min_{\ub \in W^{1,2}(\omega;\ZR^2), \ v \in W^{2,2}(\omega),\ \bar{\Rbb} \in \SO(3)} J^0 (\ub,v, \bar{\Rbb}), \textrm{ if } h^{-4}E^h \to 1,
\end{equation}
\begin{equation}
    \lim_{h \to 0} \frac{1}{E^h} J^h (\yb^h)= \min_{\ub \in W^{1,2}(\omega;\ZR^2), \ v \in W^{2,2}(\omega),\ \bar{\Rbb} \in \SO(3)} J^0_L (\ub,v, \bar{\Rbb}), \textrm{ if } h^{-4}E^h \to 0.
\end{equation}

\end{enumerate}

\end{theorem}
\begin{prooof}
The proof goes in the same direction as the proof of Theorem 2.5. in \cite{Lewicka1}.
If we take $\yb^h=\bar{\Rbb}^h x$ where $\bar{\Rbb}^h$ is chosen such that
$$ m^h=\int_{\hat{\Omega}^h} \fb^h_3(x) (\bar{\Rbb}^h x)_3 dx,  $$
we see the left inequality in (\ref{nejinf}).
Let us now take the minimizing sequence $\yb^h$. Using Theorem (\ref{prepisano}) and the coercivity property of energy density function we find $\bar{\Rbb}^h \in \SO(3)$ such that
\begin{equation}
\| (\nabla \yb^h) \circ \Thetab^h \circ P^h- \bar{\Rbb}^h\|_{L^2(\Omega;\ZR^{3 \times 3})} \leq Ch^{-1} \sqrt{I^h (\yb^h)}.
\end{equation}
Taking $\cb^h= \int_{\Omega} \yb^h \circ \Thetab^h \circ P^h $ we conclude, by using Poincare inequality the boundedness of $\nabla \Thetab^h$, $(\nabla \Thetab^h)^{-1}$, $\nabla P^h$  that
\begin{eqnarray}
\nonumber & &\| (\bar{\Rbb}^h)^T \yb^h -\cb^h -\ide \|_{L^2(\hat{\Omega}^h;\ZR^3)} \leq \\ \nonumber & & \hspace{15ex} \leq C \frac{1}{\sqrt{h}} \| (\bar{\Rbb}^h)^T (\yb^h \circ \Thetab^h \circ P^h) -\cb^h - \Thetab^h \circ P^h \|_{L^2(\Omega;\ZR^3)} \\ \nonumber & & \hspace{15ex} \leq C \frac{1}{\sqrt{h}} \| (\bar{\Rbb}^h)^T \nabla (\yb^h \circ \Thetab^h \circ P^h) -\nabla (\Thetab^h \circ P^h) \|_{L^2(\Omega;\ZR^{3 \times 3})} \\ \nonumber
& & \hspace{15ex} \leq C \frac{1}{\sqrt{h}} \| (\bar{\Rbb}^h)^T (\nabla \yb^h) \circ \Thetab^h \circ P^h -\Ibb \|_{L^2(\Omega;\ZR^{3 \times 3})} \\ \label{zzzzzz}& & \hspace{15ex} \leq C\frac{1}{\sqrt{h}} h^{-1} \sqrt{I^h (\yb^h)}.
\end{eqnarray}
Using (\ref{pretpostavkesile}) we conclude
\begin{eqnarray} \nonumber
I^h (\yb^h) &=& J^h(\yb^h)+ \frac{1}{h} \int_{\hat{\Omega}^h} \fb_3^h \yb_3^h -\frac{m^h}{h} \\ \nonumber
&=& J^h(\yb^h)+\frac{1}{h} \int_{\hat{\Omega}^h} (\bar{\Rbb}^h)^T \left( \begin{array}{c} 0 \\ 0 \\ \fb_3^h (x) \end{array} \right) \cdot ((\bar{\Rbb}^h)^T \yb^h (x) -\cb^h -x)dx\\  \label{izrazcicak} & &+ \frac{1}{h} \int_{\hat{\Omega}^h} (\bar{\Rbb}^h)^T  \left( \begin{array}{c} 0 \\ 0 \\ \fb_3^h (x) \end{array} \right) \cdot x dx -\frac{m^h}{h}.
\end{eqnarray}
\end{prooof}
Using  (\ref{pretpostavkesile}), (\ref{svojmin}), the left inequality in (\ref{nejinf}) and the definition of $m^h$  we conclude that there exists  $C>0$ such that
\begin{equation}
I^h (\yb^h)  \leq C (E^h+ \sqrt{E^h} \sqrt{I^h (\yb^h)}).
\end{equation}
From this we conclude that there exists $C>0$ such that $I^h(\yb^h) \leq CE^h$.
Using (\ref{zzzzzz}) we conclude $\frac{1}{E^h} J^h (\yb^h) \geq -C$. From (\ref{svojmin}) we conclude the right inequality in (\ref{nejinf}).

Everything else is the consequence of the fact that
$\frac{1}{E^h}J^h \stackrel{\Gamma}{\to} J^0$ (see
\cite{Braides,dalmaso}).

To prove the lower bound we take $\yb^h \in W^{1,2}(\hat{\Omega}^h;\ZR^3)$ such that \\ $\liminf_{h \to 0} \frac{1}{E^h} J^h(\yb^h) <+\infty$.
We take subsequence (not relabeled) such that $\lim_{h \to 0} \frac{1}{E^h} J^h(\yb^h) =\liminf_{h \to 0} \frac{1}{E^h} J^h(\yb^h)$. In the same way as in (\ref{izrazcicak}) we conclude $I^h(\yb^h)<CE^h$.

Using Lemma \ref{osnovna} we find $\bar{\Rbb}^h \in \SO(3)$, $\cb^h \in \ZR$, $\ub \in W^{1,2}(\Omega;\ZR^2)$, $v \in W^{2,2}(\omega)$ such that for the subsequence (not relabeled) $ \widetilde{\yb}^h:=(\bar{\Rbb}^h)^T \yb^h- \cb^h$, $\ub^h$, $v^h$  it is valid $\ub^h \rightharpoonup \ub$ in $W^{1,2}(\omega;\ZR^2)$ and $v^h \to v$ in $W^{1,2}(\omega)$. Let us take subsequence (also not relabeled) such that $\bar{\Rbb}^h \to \bar{\Rbb}$.
Let us write
\begin{eqnarray} \nonumber
\frac{1}{E^h}J^h (\yb^h)&=&\frac{1}{E^h} I^h(\widetilde{\yb}^h)+\frac{1}{hE^h}\int_{\hat{\Omega}^h} \fb_3^h (x) (\widetilde{\yb}_3^h(x)-x_3) dx\\ & &+\frac{1}{hE^h} \Big( m^h-\int_{\hat{\Omega}^h} \fb^h_3(x) (\bar{\Rbb}^h x)_3 dx\Big). \\
\end{eqnarray}
Let us prove that
\begin{equation} \label{malopoduze}
\frac{1}{hE^h} \int_{\hat{\Omega}^h} \fb_3^h (x) (\widetilde{\yb}_3^h(x)-x_3) dx \to \int_{\omega} \Big(\int_{-1/2}^{1/2} \fb_3(x',x_3)dx_3\Big) v(x') dx'.
\end{equation}
Let us consider
\begin{eqnarray} \nonumber
& &\Big|\frac{1}{hE^h}\int_{\hat{\Omega}^h} \fb_3^h (x) (\widetilde{\yb}_3^h(x)-x) dx- \int_{\omega} \Big( \int_{-1/2}^{1/2} \fb_3(x',x_3)dx_3\Big) v(x') dx' \Big| \\ \nonumber & & \leq \Big|\int_{\Omega}\frac {(\fb_3^h \circ \Thetab^h \circ P^h)(x)}{h\sqrt{E^h}} \frac{h}{\sqrt{E^h}} \Big[ ((\widetilde{\yb}_3^h \circ \Thetab^h \circ P^h)(x)-(\Thetab^h_3 \circ P^h)(x))\\ \nonumber & & \hspace{25ex} -v(x')\Big] \det\left((\nabla \Thetab^h)(P^h(x))\right)dx\Big| \\ \nonumber & &
+\Big| \int_{\Omega}\frac {(\fb_3^h \circ \Thetab^h \circ P^h)(x)}{h\sqrt{E^h}} v(x')\det\left((\nabla \Thetab^h)(P^h(x))\right)dx -\\  \label{josdugi} & &-\int_{\Omega} \fb_3(x) v(x') \det\left((\nabla \Thetab^h)(P^h(x))\right)     \Big|
+C(f^h)^2.
\end{eqnarray}
The first term in  (\ref{josdugi}) goes to $0$ by Lemma \ref{naknadnododano}, the boundedness of  $\nabla \Thetab^h$ and the weak convergence in (\ref{pretpostavkesile})
and the second term goes to $0$ by the weak convergence in (\ref{pretpostavkesile}).
Thus we have proved (\ref{malopoduze}). From the fact that $|\frac{1}{E^h}J^h (\yb^h)|$ is bounded we conclude that $\frac{1}{hE^h} \Big( m^h-\int_{\hat{\Omega}^h} \fb^h_3(x) (\Rbb^h x)_3 dx\Big)$ is bounded. From that we conclude $\bar{\Rbb} \in \mathcal{M}$.
Using Corollary (\ref{zajakuisto}) we conclude
 \begin{eqnarray*}
 \liminf_{h \to 0}\frac{1}{E^h}J^h (\yb^h) &\geq& \liminf_{h \to 0} \frac{1}{E^h} I^h(\yb^h) +\lim_{h \to 0} \frac{1}{hE^h} \int_{\hat{\Omega}^h} \fb_3^h (x) (\widetilde{\yb}_3^h(x)-x_3) dx\\ & &+\liminf_{h \to 0} \frac{1}{hE^h} \Big( m^h-\int_{\hat{\Omega}^h} \fb^h_3(x) (\bar{\Rbb}^h x)_3 dx \Big) \\
&\geq& \left\{ \begin{array}{ll} J^0 (\ub,v,\bar{\Rbb}),& \textrm{ if } h^{-4}E^h \to 1 \\ J^0_L (\ub,v,\bar{\Rbb}),& \textrm{ if } h^{-4}E^h \to 0 \end{array} \right..
\end{eqnarray*}
To prove the upper bound for $\ub \in W^{1,2}(\omega;\ZR^2)$, $v \in W^{2,2}(\omega;\ZR)$, $\bar{\Rbb} \in \mathcal{M}$ we choose $\hat{\yb}^h$ from Theorem \ref{upperbound} and $\bar{\Rbb}^h \to \bar{\Rbb}$ such that
\begin{equation} \label{orrrrr} r(\bar{\Rbb})= \lim_{h \to 0} \frac{1}{hE^h} \Big( m^h-\int_{\hat{\Omega}^h} \fb^h_3(x) (\bar{\Rbb}^h x)_3 dx\Big). \end{equation}
We define $\yb^h=\bar{\Rbb}^h \hat{\yb}^h$.
We decompose the total energy functional
\begin{eqnarray} \nonumber
\frac{1}{E^h}J^h (\yb^h)&=&\frac{1}{E^h} I^h(\hat{\yb}^h)+\frac{1}{hE^h}\int_{\hat{\Omega}^h} (\bar{\Rbb}^h)^T \left( \begin{array}{c} 0 \\ 0 \\ \fb_3^h (x) \end{array} \right)\cdot (\hat{\yb}^h(x)-x) dx\\ & &+\frac{1}{hE^h} \Big( m^h-\int_{\hat{\Omega}^h} \fb^h_3(x) (\bar{\Rbb}^h x)_3 dx\Big).  \label{dekompozicija}
\end{eqnarray}
It is easily seen, by using Lemma \ref{naknadnododano}, that
\begin{eqnarray} \nonumber
& &\frac{1}{hE^h}\int_{\hat{\Omega}^h} (\bar{\Rbb}^h)^T \left( \begin{array}{c} 0 \\ 0 \\ \fb_3^h (x) \end{array} \right)\cdot  (\hat{\yb}^h(x)-x) dx =\\ & &\nonumber \int_{\Omega}
(\bar{\Rbb}^h)^T \left( \begin{array}{c} 0 \\ 0 \\  \frac {(\fb_3^h \circ \Thetab^h \circ P^h)(x)}{h\sqrt{E^h}} \end{array} \right) \cdot
 \frac{h}{\sqrt{E^h}}  \Big((\hat{\yb}^h \circ \Thetab^h \circ P^h)(x)-\\& & \hspace{15ex}-(\Thetab^h \circ P^h)(x)\Big)\det\left((\nabla \Thetab^h)(P^h(x))\right) dx \\ \label{zzz} & & \to \int_{\Omega} (\bar{\Rbb})^T \left( \begin{array}{c} 0 \\ 0 \\  \fb_3 \end{array} \right) \cdot  \left( \begin{array}{c} 0 \\ 0 \\  v \end{array} \right)dx=\bar{\Rbb}_{33} \int_{\omega} \Big(\int_{-1/2}^{1/2} \fb_3 dx_3 \Big)  v(x')dx'.
\end{eqnarray}
From the decomposition (\ref{dekompozicija}) we conclude, by using (\ref{orrrrr}) and (\ref{zzz}),  that $\frac{1}{E^h} J^h (\yb^h) \to J^0 (\ub,v, \bar{\Rbb})$, if $h^{-4}E^h \to 1$ and  $\frac{1}{E^h} J^h (\yb^h) \to J^0_L (\ub,v, \bar{\Rbb})$, if $h^{-4}E^h \to 0$.
Thus we have proved the $\Gamma$-convergence result. Existence of the minimizer follows from the $\Gamma$-convergence theory, but can also be proven independently (see Lemma \ref{egzistencija} and the references in \cite[Section 5.12]{Ciarlet0}).
\begin{remark} \label{inplane1}
The analysis in the situation when we have also in-plane forces
(dead loads) is more complicated. The problems that appear in
these situations are similar to the ones observed for the
classical plate in \cite{Muller5}. Namely, if we want to add
in-plane forces, they should be of order $h^2$ in all the
situations $h^{-4}E^h \to \infty$ or $h^{-4}E^h \to 1$, since
$\ub$ is of order $\frac{E^h}{h^2}$. From the expression
(\ref{izrazcicak}) we can only conclude that $I^h(\yb^h) \leq
Ch^2$. Thus we do not know which model is appropriate for these
forces and we would have to impose some stability condition
(see \cite{Muller5}). The situation is better in the case
$h^{-4}E^h \to 0$ (the order of the in-plane forces should then
be $\sqrt{E^h}$), despite the fact that we can not conclude
directly  the right order of the internal energy from the
expression (\ref{izrazcicak}). This is because the body can
perform large rotation and mix the in-plane and the normal
forces. This can be prevented by imposing appropriate boundary
conditions (see \cite{Muller5}).
\end{remark}
\begin{lemma} \label{egzistencija}
For $\fb_3 \in L^2 (\Omega)$ such that $\int_{\Omega} \fb_3=0$
and  $g:\SO(3)\to \ZR \bigcup \{+\infty\}$ bounded from below,
lower semi-continuous and not identically equal to
$\{+\infty\}$ on the set $\mathcal{M}'_c$, the functionals
\begin{eqnarray}
J^0(\ub,v,\bar{\Rbb}) &=& I^{MvK}(\ub,v)- \bar{\Rbb}_{33} \int_{\omega}  \fb_3(x')   v(x')dx'+g(\bar{\Rbb}), \\
J^0_L(\ub,v,\bar{\Rbb}) &=& I^{LMvK}(\ub,v)-\bar{\Rbb}_{33} \int_{\omega}  \fb_3 (x')  v(x')dx'+g(\bar{\Rbb})
\end{eqnarray}
have the minimizers in the space $W^{1,2}(\omega;\ZR^2) \times W^{2,2}(\omega) \times \mathcal{M}_c'$, where
$\mathcal{M}'_c$ is any closed subset of the set
\begin{equation} \label{uvjetprostora}
\mathcal{M}'= \{ \Qbb \in \SO(3) | \int_{\omega} \fb_3(x')  (\bar{\Qbb}\Fbb\left( \begin{array}{c} x' \\ 0 \end{array} \right))_3 dx'=0, \quad \forall \Fbb \in \so(3) \}.
\end{equation}
\end{lemma}
\begin{prooof}
Sequentially weakly lower semi-continuity of the functionals is the direct consequence of the convexity of the form $Q_2$ and Rellich-Kondrachov compactness embedding theorems (see e.g. \cite{Adams}). The boundedness of the minimizing sequence is the only fact that has to be proved. We shall do it for the functional $J^0$.
Let us take the minimizing sequence $(\ub^h, v^h,\bar{\Rbb}^h)$ such that $J^0(\ub^h,v^h,\bar{\Rbb}^h)\to \inf J^0 <+\infty$. Since
$Q_2$ is positive definite on symmetric matrices,
 we have that $I^{MvK}(\ub,v) \geq C_1 \| (\nabla')^2 v\|_{L^2(\Omega)}$. By the Poincare inequality there exist $\ab^h \in \ZR^2,\bbb^h \in \ZR$ such that $\|v^h -\ab^h \cdot x'-\bbb^h\|_{W^{2,2}(\Omega)} \leq C_2 \| (\nabla')^2 v^h \|_{L^2(\Omega)}$. Let us note that the condition in the definition of the space (\ref{uvjetprostora}) for $\Qbb_{33} \neq 0$ implies $\int_{\omega} \fb_3 x_1=0$ and $\int_{\omega} \fb_3 x_1=0$. Thus in either case we have $\bar{\Rbb}^h_{33} \int_{\omega}  \fb_3(x') (\ab^h\cdot x'+\bbb^h)dx'=0$.  From this we can easily deduce (by the standard argumentation) that $\|v^h-\ab^h \cdot x'-\bbb^h\|_{W^{2,2}(\Omega)} \leq C$. The problem is, however, that the functions of the type $\ab^h \cdot x'$ affect the part of the functional
$Q_2(\frac{1}{2}(\nabla' \ub^h+(\nabla' \ub^h)^T+\nabla' v^h \otimes \nabla' v^h+\nabla' v^h \otimes \nabla' \theta+\nabla' \theta \otimes \nabla' v^h))$. But we have
\begin{eqnarray*}
& &\nabla' \ub^h+(\nabla' \ub^h)^T+\nabla' v^h \otimes \nabla' v^h+\nabla' v^h \otimes \nabla' \theta+\nabla' \theta \otimes \nabla' v^h= \\ & &\nabla' \ub_c^h+\nabla' (\ub_c^h )^T+\nabla' v^h_c\otimes \nabla' v^h_c
+\nabla' v^h_c \otimes \nabla' \theta+\nabla' \theta \otimes \nabla' v^h_c,
\end{eqnarray*}
where
$$ \ub^h_c=\ub^h +(v^h+\theta) \ab^h+ \frac{1}{2} (\ab^h \cdot x') \ab^h ,\quad v^h_c= v^h-\ab^h \cdot x'-\bbb^h.$$
Thus we conclude that
$$J^0(\ub_c^h,v^h_c,\bar{\Rbb}^h)=J^0(\ub^h,v^h,\bar{\Rbb}^h) \to \inf J^0,$$
for $\|v^h_c\|_{W^{2,2}(\Omega)} \leq C$.
Now we take $\Abb^h_c \in \ZR^{2 \times 2}$, $\bbb^h_c\in \ZR$ such that
we have by Korn's inequality
\begin{equation} \label{Kornovapon} \|\ub_c^h-\Abb^h_c x'-\bbb^h_c\|_{W^{1,2}(\Omega;\ZR^2)} \leq C_3  \|\sym (\nabla'\ub_c^h) \|_{L^{2}(\Omega;\ZR^2)}. \end{equation}
If we define $\ub^h_{cc}=\ub_c^h-\Abb^h_c x'-\bbb^h_c$ we have
$$ J^0(\ub_{cc}^h,v^h_c,\bar{\Rbb}^h)=J^0(\ub^h,v^h,\bar{\Rbb}^h) \to \inf J^0,$$
for $\|v^h_c\|_{W^{2,2}(\Omega)} \leq C$.
From the boundedness of $v^h_c$ and the fact that
\begin{eqnarray*} I^{MvK}(\ub^h_{cc},v^h_c,\bar{\Rbb}^h) &\geq& C_1 \|\nabla' \ub_{cc}^h+\nabla' (\ub_{cc}^h )^T+\nabla' v^h_c\otimes \nabla' v^h_c+\\ & & \hspace{15ex}
+\nabla' v^h_c \otimes \nabla' \theta+\nabla' \theta \otimes \nabla' v^h_c\|_{L^2(\Omega;\ZR^2)},  \end{eqnarray*}
we easily conclude the boundedness of
$ \|\nabla' \ub_{cc}^h+\nabla' (\ub_{cc}^h )^T \|_{L^2(\Omega;\ZR^2)}$ and from (\ref{Kornovapon}) the boundedness of
$\|\ub_{cc}^h\|_{W^{1,2}(\Omega;\ZR^2)}$. From this it follows that the sequence $(\ub_{cc}^h,\vb^h,\bar{\Rbb}^h)$ has its weak limit. Due to the sequentially weakly lower semi-continuity of the functional $J^0$, the weak limit of the minimizing sequence is a minimizer.
\end{prooof}
\begin{remark} \label{inplane2}
If we want to include tangential forces and displacement
boundary conditions,  then we could obtain existence result for
the functional $J^0$ only for small tangential forces (see
\cite{Ciarlet2,Gratie,Ciarlet3}). This is not the case for the
functional $J^0_L$ where we can, by using the generalized
Korn's inequality (see \cite{Ciarlet0}), conclude the
coercivity of the functional in the appropriate space including
the displacement boundary conditions.
\end{remark}


\begin{thebibliography}{XX}
\bibitem{Ciarlet0}  Ciarlet, P.G.: Mathematical elasticity.
    Vol. II. Theory of plates. Studies in Mathematics and its
    Applications, 27. North-Holland Publishing Co., Amsterdam
    (1997).
\bibitem{Vorovich}  Vorovich, I.I: Nonlinear Theory of Shallow
    Shells, Springer-Verlag: New York, 1999.

\bibitem{Ciarlet1}  Ciarlet, P.G.: Mathematical elasticity.
    Vol. III. Theory of shells. Studies in Mathematics and its
    Applications, 29. North-Holland Publishing Co., Amsterdam,
    2000.

\bibitem{Andreoiu}   Andreoiu, G., Faou, E.: Complete
    asymptotics for shallow shells, Asymptotic Analysis 25 ,
    239--270 (2001).
\bibitem{Fox}  Fox, D.D., Raoult A.,  Simo, J.C.: A
    justification of nonlinear properly invariant plate
    theories, Arch. Rational Mech. Anal., 124, p. 157–199
    (1993).
\bibitem{Andreoiu2}  Andreoiu-Banica, G.: Justification of the
    Marguerre-von K\'{a}rm\'{a}n equations in curvilinear
    coordinates, Asymptotic Analysis 19, 35--55 (1999).

\bibitem{Hamdouni} Hamdouni, A., Millet, O.: Classification of
    thin shell models deduced from the nonlinear
    three-dimensional elasticity. Part I : the shallow shells,
    Arch. Mech., 55, 135--175 (2003).
\bibitem{Ciarlet2} Ciarlet, P.G., Gratie, L.: On the existence
    of solutions to the generalized Marguerre-von
    K\'{a}rm\'{a}n equations, Math. Mech. Solids  11, 83--100
    (2006).
\bibitem{Gratie} Gratie, L.: Generalized Marguerre-von
    K\'{a}rm\'{a}n equations of a nonlinearly elastic shallow
    shell, Applicable Analysis, 81(5), 1107-1126 (2002).
\bibitem{Braides} A. Braides: $\Gamma$-convergence for
    Beginners, Oxford University Press, Oxford, 2002.
\bibitem{dalmaso} Dal Maso, G.:  An introduction to
    $\Gamma$-convergence, Progress in Nonlinear Differential
    Equations and Their Applications, Birk\"auser, Basel, 1993.
\bibitem{Ledret1} Le Dret, H.,  Raoult, A.: The nonlinear
    membrane model as a variational limit of nonlinear
    three-dimensional elasticity, Journal de Math\'ematiques
    Pures et Appliqu\'ees 74, 549--578 (1995).

\bibitem{LeDret2} Le Dret, H.,  Raoult, A.: The membrane shell
    model in nonlinear elasticity: A variational asymptotic
    derivation, Journal of Nonlinear Science 6, Number 1,
    59--84 (1996).
\bibitem{Muller0}  Friesecke, G.,  James R.D.,  M\"uler, S.: A
    theorem on geometric rigidity and the derivation of
    nonlinear plate theory from three-dimensional elasticity,
    Comm. Pure Appl. Math.  55, 1461--1506 (2002).
\bibitem{Muller3}  Friesecke, G., James R.D.,  M\"uler, S.: A
    Hierarchy of Plate Models Derived from Nonlinear Elasticity
    by $\Gamma$-Convergence, Archive for Rational Mechanics and
    Analysis  180, no.2, 183--236 (2006).
\bibitem{Muller4}  Friesecke, G.,  James R.D.,  M\"uler, S.:
    The F\"oppl-von K\'{a}rm\'{a}n plate theory as a low energy
    $\Gamma$-limit of nonlinear elasticity, Comptes Rendus
    Mathematique 335, no. 2, 201--206  (2002).

\bibitem{Muller6} Friesecke, G., James R., Mora, M.G.,
    M\"uller, S.: Derivation of nonlinear bending theory for
    shells from three-dimensional nonlinear elasticity by
    $\Gamma$-convergence, C. R. Math. Acad. Sci. Paris, 336,
    no. 8, 697–702 (2003).
\bibitem{Lewicka1}  Lewicka, M.,  Mora, M.G., Pakzad, M.: Shell
    theories arising as low energy $\Gamma$-limit of 3d
    nonlinear elasticity , Ann. Scuola Norm. Sup. Pisa Cl. Sci.
    5, Vol. IX, 1--43 (2010).
\bibitem{Muller5}  Lecumberry, M.,  M\"uller, S.: Stability of
    slender bodies under compression and validity of the von
    K\'{a}rm\'{a}n theory, Archive for Rational Mechanics and
    Analysis, Volume 193, Number 2, 255--310 (2009).
\bibitem{Mora}  Mora, M.G., Scardia, L.: Convergence of
    equilibria of thin elastic plates under physical growth
    conditions for the energy density, submitted paper.
\bibitem{Muller7}  M\"uller, S.,  Packzad, M.R.: Convergence of
    equilibria of thin elastic plates : the von K\'{a}rm\'{a}n
    case, Communications in Partial Differential Equations 33,
    Number 6, 1018--1032 (2008).
\bibitem{Adams} Adams, R.A.: Sobolev spaces, Academic press,
    New York 1975.
\bibitem{Lewicka2} Lewicka, M.,  Mora, M.G., Pakzad, M.: The
    matching property of infinitesimal isometries on elliptic
    surfaces and elasticity of thin shells, accepted in Arch.
    Rational Mech. Anal.


\bibitem{Lewicka3}  Lewicka, M.,  Pakzad, M.: The infinite
    hierarchy of elastic shell models: some recent results and
    a conjecture , accepted in Fields Institute Communications
    (2010).


\bibitem{Ciarlet3} Ciarlet, P.G., Gratie, L.: From the
    classical to the generalized von K\'{a}rm\'{a}n and
    Marguerre-von K\'{a}rm\'{a}n equations.  J. Comput. Appl.
    Math.  190,  no. 1-2, 470--486 (2006).

\end{thebibliography}
\end{document}